\documentclass{article}
\usepackage[utf8]{inputenc}
\usepackage{amssymb,amsmath,amsthm}
\usepackage[square]{natbib}
\usepackage{tikz,pgfplots}

\usepackage[breaklinks]{hyperref}
\usepackage{xcolor}
\hypersetup{
  colorlinks,
  linkcolor={red!50!black},
  citecolor={blue!50!black},
  urlcolor={blue!80!black}}

\allowdisplaybreaks

\newtheorem{thm}{Theorem}
\newtheorem*{thm*}{Theorem}
\newtheorem{prop}{Proposition}
\newtheorem{lem}{Lemma}
\newtheorem{cor}{Corollary}
\newtheorem{defn}{Definition}
\newtheorem{example}[thm]{Example}
\newtheorem{rem}{Remark}
\numberwithin{equation}{section}

\newcommand{\R}{\mathbb{R}}
\newcommand{\scal}[2]{\left\langle #1,#2 \right\rangle}
\DeclareMathOperator{\Anti}{Anti}

\title{Rank optimality for the Burer-Monteiro factorization}
\author{Ir\`ene Waldspurger\thanks{CNRS, Universit\'e Paris Dauphine, Inria Mokaplan, France
  (\texttt{waldspurger@ceremade.dauphine.fr}).}
\and Alden Waters\thanks{Bernoulli Institute, Rijksuniversiteit Groningen, Groningen, Netherlands (\texttt{a.m.s.waters@rug.nl}).}}

\begin{document}

\maketitle

\begin{abstract}
  When solving large scale semidefinite programs that admit a low-rank solution, an efficient heuristic is the Burer-Monteiro factorization: instead of optimizing over the full matrix, one optimizes over its low-rank factors. This reduces the number of variables to optimize, but destroys the convexity of the problem, thus possibly introducing spurious second-order critical points. The article \citep*{boumal_voroninski_bandeira} shows that when the size of the factors is of the order of the square root of the number of linear constraints, this does not happen: for almost any cost matrix, second-order critical points are global solutions. In this article, we show that this result is essentially tight: for smaller values of the size, second-order critical points are not generically optimal, even when the global solution is rank $1$.
\end{abstract}

\section{Introduction}

We consider a semidefinite program:
\begin{align}
  \mbox{minimize }&\mathrm{Trace}(CX)\tag{SDP}\label{eq:SDP_intro}\\
  \mbox{such that }&\mathcal{A}(X)=b,\nonumber\\
                  &X\succeq 0,\nonumber
\end{align}
where the variable $X$ and the fixed matrix $C$ are symmetric, of size $n\times n$, and $\mathcal{A}$ is a linear operator capturing $m$ equality constraints.

Various iterative algorithms have been developed to solve such a problem at a given precision level, but tend to be computationally demanding. For example, in full generality, each iteration may cost $O((m+n)mn^2)$ arithmetic operations with an interior-point solver \citep*[Page 357]{borchers} (assuming $m\leq \mathrm{dim}(\mathrm{Sym}(n)) = \frac{n(n+1)}{2}$), and $O((m+n)n^2)$ with first-order techniques applied to a smoothed version of the problem \citep*[Section 3]{nesterov05}.

Improvements are possible if $\mathcal{A}$ has some structure that can be exploited, but they often do not suffice to make large-scale semidefinite programs computationally easy. Another property can then be used: semidefinite programs tend to have a low-rank minimizer (in many applications, there is one with rank $O(1)$, and, in any case, always one with rank $\sim \sqrt{2m}$ \citep*[Theorem 2.1]{pataki}). Low-rank matrices can be stored and manipulated in a much more efficient way than full-rank ones, which allows for less computationally demanding algorithms.

Frank-Wolfe methods, in particular, take advantage of this \citep*{jaggi, laue, yurtsever}. Here, we are interested in another approach, the Burer-Monteiro factorization \citep*{burer}. Its principle is that a semidefinite matrix with rank $p\ll n$ can be factorized as $
X = U U^T, $
with $U\in\R^{n\times p}$. Assuming that a low-rank solution $X_{opt}$ exists, if $p\geq \mathrm{rank}(X_{opt})$, Problem \eqref{eq:SDP_intro} is then equivalent to
\begin{align}
  \mbox{minimize }&\mathrm{Trace}(CUU^T)\tag{Factorized SDP}\label{eq:Factorized_SDP_intro}\\
  \mbox{such that }&\mathcal{A}(UU^T)=b,\nonumber\\
                  &U\in\R^{n\times p}.\nonumber
\end{align}
Now the unknown $U$ has $np$ coordinates, fewer than the $n^2$ coordinates of $X$. Consequently, we can run on Problem \eqref{eq:Factorized_SDP_intro} local optimization algorithms that would be too slow on Problem \eqref{eq:SDP_intro}. The caveat is that, since the factorized problem is not convex, they are not guaranteed to find a global minimizer, at best a second-order critical point. Nevertheless, they work extremely well in many applications. Typically, as soon as $p$ is slightly larger than $\mathrm{rank}(X_{opt})$, local optimization algorithms seem to globally solve Problem \eqref{eq:P}. Numerical examples where this phenomenon occurs can be found in \citep*{burer_monteiro03}, \citep*[Section 5]{journee}, \citep*[Section 5]{boumal_riemannian} or \citep*[Section 5]{rosen_carlone_bandeira_leonard}.
% Rq : les figures 1 et 2 de [Journée et al] montrent qu'on peut prendre un rang assez petit pour les problèmes MaxCut étudiés dans [Burer-Monteiro 2003], bien plus petit que celui utilisé par Burer et Monteiro.

The article \citep*{bandeira_low_rank} gives a rigorous explanation of this behavior for instances of \eqref{eq:Factorized_SDP_intro} coming from $\mathbb{Z}_2$-synchronization and community detection. It notably establishes, in particular statistical regimes where Problem \eqref{eq:SDP_intro} has a rank-$1$ solution, that all second-order critical points of Problem \eqref{eq:Factorized_SDP_intro} with $p=2$ are global minimizers. Hence, suitable local optimization algorithms globally solve Problem \eqref{eq:Factorized_SDP_intro}. Similarly, \citep*{ge_lee_ma,sun_qu_wright,li_zhu_tang} show, in related settings, that all second-order critical points of the Burer-Monteiro factorization are the optimal solution as soon as $p\geq \mathrm{rank}(X_{opt})$.

But these works apply in very specific settings only. They provide no general theory on when local optimization algorithms solve Problem \eqref{eq:Factorized_SDP_intro}. With no restrictive assumptions, essentially the only result is \citep*{boumal_voroninski_bandeira}. Building on \citep*{burer} and \citep*{boumal_riemannian}, it shows that, under simple hypotheses, all second-order critical points of Problem \eqref{eq:Factorized_SDP_intro} are global minimizers, for almost any matrix $C$, as soon as
\begin{equation}\label{eq:lower_bound_p_intro}
  \frac{p(p+1)}{2} >m,
\end{equation}
that is $p> \lfloor\sqrt{2m + 1/4}-1/2\rfloor$. Extensions can be found in \citep*{pumir_jelassi_boumal,  bhojanapalli_boumal_jain_netrapalli}.

As a result there is a gap in the literature: in all the concrete settings that could be studied, all second-order critical points of Problem \eqref{eq:Factorized_SDP_intro} are global minimizers as soon as $p\gtrsim \mathrm{rank}(X_{opt})$, in line with numerical experiments, but in the general case, the only guarantees at our disposal state that we need $p$ to be at least as large as $\sim \sqrt{2m}$. In many applications, $\mathrm{rank}(X_{opt})=O(1)$ while $m=O(n)$, hence these two estimates are far apart, making a huge difference on the computational cost of certifiable algorithms.

The natural question is, "can the gap be reduced?" In this article, we negatively answer this question, and show that Inequality \eqref{eq:lower_bound_p_intro} is essentially optimal.
It can be slightly improved, to
%Theorem \ref{theorem:bvb_improved} shows a minor improvement over Inequality \eqref{eq:lower_bound_p_intro} is possible: Under a stronger (but reasonable) geometrical assumption than in \citep*{boumal_voroninski_bandeira}, all second-order critical points of Problem \eqref{eq:Factorized_SDP_intro} are global minimizers, for almost any $C$, as soon as
  \begin{equation}\label{eq:lower_bound_improved_intro}
    \frac{p(p+1)}{2}+p >m.
  \end{equation}
  %under a stronger (but reasonable) geometrical assumption than in \citep*{boumal_voroninski_bandeira}.
  This is Theorem \ref{theorem:bvb_improved}.
  But Theorem \ref{theorem:main} (our main result) states that, under reasonable assumptions on $\mathcal{A},b$, if $p$ is such that
  \begin{equation*}
    \frac{p(p+1)}{2} + pr_* \leq m,
  \end{equation*}
  where $r_*=\min\{\mathrm{rank}(X),X\succeq 0,\mathcal{A}(X)=b\}$, there exists a set of cost matrices $C$ with non-zero Lebesgue measure on which Problem \eqref{eq:SDP_intro} admits a global minimizer with rank $r_*$, but Problem \eqref{eq:Factorized_SDP_intro} has second-order critical points which are not global minimizers. In particular, if $r_*=1$ (as is the case in \textit{MaxCut} relaxations, for instance), Inequality \eqref{eq:lower_bound_improved_intro} is exactly optimal. Therefore without specific assumptions on $C$, when running a local optimization algorithm on Problem \eqref{eq:Factorized_SDP_intro} with $p$ smaller than $\sim\sqrt{2m}$, we cannot be sure not to run into a spurious second-order critical point, even if there exists a global minimizer with rank $O(1)$.

  Regarding the organization of this article, Section \ref{s:preliminaries} contains basic definitions (Subsection \ref{ss:definitions}) and properties (Subsection \ref{ss:basic}), defines and discusses an important assumption for our main result (Subsection \ref{ss:min_sec}). Section \ref{s:main_results} presents the main results: Theorems \ref{theorem:bvb_improved} and \ref{theorem:main} are respectively stated in Subsections \ref{ss:bvb_improved} and \ref{ss:main}. Subsection \ref{ss:examples} provides examples. The other sections contain the proofs: Theorem \ref{theorem:bvb_improved} is proved in Section \ref{s:proof_bvb_improved}, and Theorem \ref{theorem:main} in Section \ref{s:proof}. 

\subsection{Notation}

For any $p,q\in\mathbb{N}^*$, we denote by $I_p$ the $p\times p$ identity matrix, and by $0_{p, q}$ the zero $p\times q$ matrix. For any $p\in\mathbb{N}^*$, we denote by $\mathbb{S}^{p\times p}$ the set of real symmetric $p\times p$ matrices, by $\Anti(p)$ the set of antisymmetric $p\times p$ matrices, and by $O(p)$ the set of orthogonal $p\times p$ matrices. For any $n_1,n_2$, we equip $\R^{n_1\times n_2}$, the set of $n_1\times n_2$ matrices, with the usual scalar product:
\begin{equation*}
  \forall M_1,M_2\in\R^{n_1\times n_2},\quad
  \scal{M_1}{M_2}\overset{def}{=} \mathrm{Tr}(M_1^T M_2).
\end{equation*}
The same formula also defines a scalar product on $\mathbb{S}^{p\times p}$, for any $p\in\mathbb{N}^*$. In both cases, the associated norm is the Frobenius norm, which we denote by $||.||_F$. For any $p\in\mathbb{N}^*$, we define $\mathrm{diag}:\R^{p\times p}\to\R^p$ as the operator which associates to a matrix the vector of its diagonal elements.

For any element $x$ of a metric space, and any positive $\epsilon$, we denote $B(x,\epsilon)$ the open ball with radius $\epsilon$, and $\overline{B}(x,\epsilon)$ the closed ball.
When $\mathcal{M}$ is a manifold, and $x$ an element of $\mathcal{M}$, we denote by $T_x\mathcal{M}$ the tangent space of $\mathcal{M}$ at $x$.

\section{Preliminaries\label{s:preliminaries}}

\subsection{Definitions\label{ss:definitions}}

We consider a problem of the following form:
\begin{align}
  \mbox{minimize }&\scal{C}{X}\label{eq:SDP}\tag{SDP}\\
  \mbox{such that }    &\mathcal{A}(X)=b,\nonumber\\
                  &X\succeq 0.\nonumber
\end{align}
Here, $\mathcal{A}:\mathbb{S}^{n\times n}\to\R^m$ is a fixed linear map, $b$ a fixed element of $\R^m$, and $C$ an element of $\mathbb{S}^{n\times n}$, which is called the \emph{cost matrix}.

We denote by $\mathcal{C}$ the set of feasible points for this problem:
\begin{equation*}
  \mathcal{C}=\{ X\in\mathbb{S}^{n\times n}, \mathcal{A}(X)=b,X\succeq 0\}.
\end{equation*}

As explained in the introduction, if we assume that Problem \eqref{eq:SDP} has an optimal solution $X_{opt}$ with rank $r$, and fix some $p\geq r$, it is equivalent to its \emph{rank $p$ Burer-Monteiro factorization}:
\begin{align}
  \mbox{minimize }&\scal{C}{VV^T}\label{eq:P}\tag{Factorized SDP}\\
  \mbox{such that }&\mathcal{A}(VV^T)=b,\nonumber\\
  &V\in\R^{n\times p}.\nonumber
\end{align}

We denote by $\mathcal{M}_p$ the set of feasible points for the factorized problem:
\begin{equation*}
  \mathcal{M}_p = \{ V\in\R^{n\times p}, \mathcal{A}(VV^T)=b\}.
\end{equation*}
It is invariant under multiplication by elements of $O(p)$. We assume that it is sufficiently regular so that we can apply smooth optimization algorithms to Problem \eqref{eq:P}. More precisely, all our results require that $(\mathcal{A},b)$ is \emph{$p$-regular}:
\begin{defn}
  For some $p\in\mathbb{N}^*$, $(\mathcal{A},b)$ is said to be $p$-\emph{regular} if, for all $V\in\mathcal{M}_p$, the following linear map is surjective:
  \begin{equation*}
    \dot{V} \in \R^{n\times p} \to \mathcal{A}(V \dot{V}^T + \dot{V}V^T) \in \R^m.
  \end{equation*}
\end{defn}
This assumption is of the same style as \citep*[Assumption 1.1]{boumal_voroninski_bandeira}. It notably guarantees \citep*[Proposition 3.3.3]{absil} that $\mathcal{M}_p$ is a submanifold of $\R^{n\times p}$, with dimension $
\dim(\mathcal{M}_p) = np - m$,
whose tangent space at any point $V$ is
\begin{equation*}
  T_V \mathcal{M}_p = \{\dot{V}\in\R^{n\times p}, \mathcal{A}(V \dot{V}^T + \dot{V}V^T) = 0\}.
\end{equation*}
The scalar product of $\R^{n\times p}$ defines a metric on the manifold $\mathcal{M}_p$, which we then view as a Riemannian manifold. Many algorithms exist for attempting to minimize a smooth function on a Riemannian manifold; a classical reference on this topic is \citep*{absil}.

However, they are a priori not guaranteed to find a global minimizer of Problem \eqref{eq:P}, but only an (approximate) first or second-order critical point of the cost function $V\in\mathcal{M}_p\to \scal{C}{VV^T}$ \citep*{boumal_absil_cartis}. These points are defined as follows:
\begin{defn}
  Let $\mathcal{N}$ be a Riemannian manifold, and $f:\mathcal{N}\to \R$ a smooth function. We denote $\nabla$ and $\mathrm{Hess}$ its gradient and Hessian with respect to the manifold.

  For any $x_0\in\mathcal{N}$, we say that $x_0$ is a \emph{first-order critical point} of $f$ if $\nabla f(x_0) = 0$  and a \emph{second-order critical point} of $f$ if
 $\nabla f(x_0) = 0$ and $\mathrm{Hess}f(x_0) \succeq 0.$
\end{defn}
The goal of this article is to study for which values of $p$ the set of second-order critical points coincides with the set of global minimizers of Problem \eqref{eq:P}.

We note that there are pairs $(\mathcal{A},b)$ which are not $p$-regular, regardless of the value of $p$ (an example is when $0_{n,n}\in\mathcal{C}$). This setting is significantly different from the one that we consider in this article: $\mathcal{M}_p$ may then have singularities, and classical Riemannian tools are a priori not applicable to Problem \eqref{eq:P}.

\subsection{Basic properties\label{ss:basic}}

It is convenient to be able to describe the solutions of Problem \eqref{eq:SDP} in terms of Karush-Kuhn-Tucker conditions. This is a priori possible only if strong duality holds but, fortunately for us, strong duality always holds when $(\mathcal{A},b)$ is $p$-regular for some $p$, yielding the following proposition (whose proof is in Appendix \ref{app:KKT_SDP}).
\begin{prop}\label{prop:KKT_SDP}
  We assume that there exists $p\in\mathbb{N}^*$ such that $(\mathcal{A},b)$ is $p$-regular and $\mathcal{M}_p\ne\emptyset$. Then a matrix $X_0\in\mathcal{C}$ is a solution of Problem \eqref{eq:SDP} if and only if there exist $g_1\in\mathbb{R}^m,C_1\in\mathbb{S}^{n\times n}$ such that
  \begin{itemize}
  \item $C=\mathcal{A}^*(g_1)+C_1$;
  \item $C_1\succeq 0$;
  \item $C_1X_0=0$.
  \end{itemize}
\end{prop}
When, $g_1,C_1$ satisfy $\mathrm{rank}(C_1)=n-\mathrm{rank}(X_0)$ in addition to the above three conditions, we say that \textit{strict complementary slackness} holds. The following proposition (whose proof is in Appendix \ref{app:strict_complementary_slackness}) states that, under an additional condition on $X_0$, it implies that the solution of Problem \eqref{eq:SDP} is unique.
\begin{prop}\label{prop:strict_complementary_slackness}
  If strict complementary slackness holds and $X_0$ is an extremal point of $\mathcal{C}$, then $X_0$ is the unique solution of Problem \eqref{eq:SDP}.
\end{prop}

The next proposition characterizes, in a similar way as Proposition \ref{prop:KKT_SDP}, the first-order critical points of Problem \eqref{eq:P}. Its proof is in Appendix \ref{app:first_order}.
\begin{prop}\label{prop:first_order}
  We assume that $(\mathcal{A},b)$ is $p$-regular for some $p\in\mathbb{N}^*$. A matrix $V\in\mathcal{M}_p$ is a first-order critical point of Problem \eqref{eq:P} if and only if there exist $g_2\in\mathbb{R}^m,C_2\in\mathbb{S}^{n\times n}$ such that
  \begin{itemize}
  \item $C=\mathcal{A}^*(g_2)+C_2$;
  \item $C_2V=0$.
  \end{itemize}
  When it exists, the pair $(g_2,C_2)$ is unique.
\end{prop}
Finally, we also provide a reformulation of second-order criticality; the proof is in Appendix \ref{app:second_order}.
\begin{prop}\label{prop:second_order}
  We assume that $(\mathcal{A},b)$ is $p$-regular for some $p\in\mathbb{N}^*$. Let $V\in\mathcal{M}_p$ be a first-order critical point of Problem \eqref{eq:P}, whose cost function we denote $f_C$. For any $\dot{V}\in T_V\mathcal{M}_p$,
  \begin{equation*}
    \mathrm{Hess} f_C(V)\cdot (\dot V,\dot V) = 2 \scal{C_2}{\dot V\dot V^T},
  \end{equation*}
  with $C_2$ defined as in Proposition \ref{prop:first_order}.
  Thus, $V$ is second-order critical if and only if
  \begin{equation}\label{eq:second_order}
    \forall\dot V\in T_V\mathcal{M}_p,\quad \scal{C_2}{\dot V\dot V^T}\geq 0,
  \end{equation}
\end{prop}
Using the notation of the previous proposition, we observe that, since $f_C$ is invariant under right multiplication by elements of $O(p)$, $\mathrm{Hess}f_C(V)\cdot(\dot V,\dot V)=0$ for any $\dot V$ tangent to the orbit of $V$ under the action of $O(p)$, that is $\dot V=VA$ for some $A\in\Anti(p)$. This motivates the following definition.
\begin{defn}
  A second-order critical point $V$ of Problem \eqref{eq:P} is \emph{non-degenerate} if, in Equation \eqref{eq:second_order}, the equality is attained exactly for matrices $\dot{V}$ of the form $\dot{V}=VA,A\in\Anti(p)$.
\end{defn}

\begin{rem}\label{rem:non_degenerate_rank}
  Equivalently, a second-order critical point is non-degenerate if
  \begin{align*}
    \mathrm{rank}(\mathrm{Hess}f_C(V))
    &= \dim(\mathcal{M}_p)-\dim\{VA,A\in\Anti(p)\} \\
    &= \dim(\mathcal{M}_p) - \frac{p(p-1)}{2}.
  \end{align*}
\end{rem}

\subsection{Definition of ``minimally secant''\label{ss:min_sec}}

The following technical property is needed for our main theorem.
\begin{defn}\label{def:min_sec}
  Let $p\in\mathbb{N}^*$ be such that $(\mathcal{A},b)$ is $p$-regular. Let $r$ be in $\mathbb{N}^*$.
  
  Let $X_0$ be a rank $r$ element of $\mathcal{C}$ and $V$ be in $\mathcal{M}_p$. We say that $\mathcal{M}_p$ is \emph{$X_0$-minimally secant at $V$} if the following three conditions hold:
  \begin{enumerate}
  \item $\mathrm{rank}(V)=p$;\label{def:min_sec1}
  \item $\mathrm{Range}(X_0) \cap \mathrm{Range}(V) = \{0\}$;\label{def:min_sec2}
  \item for any $\dot V\in T_V\mathcal{M}_p$,
    if $\mathrm{Range}(\dot V) \subset \mathrm{Range}(X_0) + \mathrm{Range}(V)$, then $\dot V = V A$ for some $A \in\Anti(p)$.\label{def:min_sec3}
  \end{enumerate}
\end{defn}
We observe that, for any $V\in\mathcal{M}_p$, the intersection
\begin{equation*}
  T_V\mathcal{M}_p
  \cap \{\dot V\in\mathbb{R}^{n\times p},\mathrm{Range}(\dot V) \subset
  \mathrm{Range}(X_0) + \mathrm{Range}(V) \}
\end{equation*}
necessarily contains $\{VA,A\in\Anti(p)\}$. Therefore, the third property in the above definition amounts to requiring that the intersection is ``as small as possible'' (hence the name ``minimally secant'').

Our main result, Theorem \ref{theorem:main}, contains the assumption that there exists $X_0,V$ such that $\mathcal{M}_p$ is $X_0$-minimally secant at $V$. This assumption is crucial for the proof: to show the existence of at least one cost matrix $C$ for which a spurious second-order critical point exists, our strategy is to fix $X_0$ and $V$, and construct $C$ for which a second-order critical point exists and is precisely $V$, while the global minimizer corresponds to $X_0$. To ensure first and second-order criticality for $V$, and global optimality for $X_0$, we essentially need that the quadratic form defined by $C$ satisfies some properties when restricted to $\mathrm{Range}(X_0)$, some other properties on $\mathrm{Range}(V)$, and still some other ones on $\mathrm{Range}(\dot V)$ for $\dot V\in T_V\mathcal{M}_p$. In order for these properties to be compatible with each other, the spaces $\mathrm{Range}(X_0),\mathrm{Range}(V),\mathrm{Range}(\dot V)$ must ``not intersect too much''. The formal content behind ``not intersecting too much'' is precisely Definition \ref{def:min_sec}.

However, as explained in Appendix \ref{app:min_sec}, when $\frac{p(p+1)}{2}+pr \leq m$, we expect such $X_0,V$ to almost always exist. This is notably the case for \textit{MaxCut} problems (see Paragraph \ref{sss:maxcut} and Appendix \ref{app:maxcut}), and also for \textit{Orthogonal-Cut} (Paragraph \ref{sss:orthogonal_cut}) and optimization over a product of spheres (Paragraph \ref{sss:product_of_spheres}).

\section{Main results\label{s:main_results}}

\subsection{Regime where critical points are global minimizers\label{ss:bvb_improved}}

As previously stated, most smooth optimization algorithms, applied to Problem \eqref{eq:P}, are only guaranteed to find a critical point of this problem, and not a global minimizer. Fortunately, \citep*{boumal_voroninski_bandeira} shows that, when $p$ is large enough, second-order critical points are always global minimizers, for almost all cost matrices $C$. Therefore, algorithms able to find second-order critical points (for instance, the trust-region method) actually solve Problem \eqref{eq:P} to optimality, provided that $C$ is ``generic''. A restated version of the theorem in \citep*{boumal_voroninski_bandeira}, under minor modifications, is the following:

\begin{thm*}\citep*[Theorem 1.4]{boumal_voroninski_bandeira}
Let $p\in\mathbb{N}^*$ be fixed. We assume that
\begin{enumerate}
\item The set $\mathcal{C}$ of feasible points for Problem \eqref{eq:SDP} is compact;
\item $(\mathcal{A},b)$ is $p$-regular.
\end{enumerate}
If
\begin{equation}\label{eq:bvb}
  \frac{p(p+1)}{2} > m,
  \qquad\left(\iff p>\left\lfloor \sqrt{2m + \frac{1}{4}} - \frac{1}{2} \right\rfloor\right)
\end{equation}
then, for almost all cost matrices $C\in \mathbb{S}^{n\times n}$, if $V\in\mathcal{M}_p$ is a second-order critical point of Problem \eqref{eq:P}, then
\begin{itemize}
\item $V$ is a global minimizer of Problem \eqref{eq:P};
\item $X=VV^T$ is a global minimizer of Problem \eqref{eq:SDP}.
\end{itemize}
\end{thm*}

It is natural to ask whether Condition \eqref{eq:bvb} is optimal, or whether the same guarantees hold for smaller ranks $p$, allowing further reductions in the computational complexity of solving Problem \eqref{eq:P}.
Our first result is that Condition \eqref{eq:bvb} can be slightly relaxed.
%, at the price of an additional assumption on $\mathcal{M}_p$.
\begin{thm}\label{theorem:bvb_improved}
  Let $p\in\mathbb{N}^*$ be fixed. We assume that
  \begin{enumerate}
  \item The set $\mathcal{C}$ of feasible points for Problem \eqref{eq:SDP} is compact.
  \item $(\mathcal{A},b)$ is $p$-regular;
  %\item If $\frac{p(p+1)}{2} \leq m$, then $\mathcal{M}_p$ is face regular.
  \end{enumerate}
  If $
    \frac{p(p+1)}{2}+p > m,$
  then the same conclusion holds as in the previous theorem.
\end{thm}

% Compared to the theorem in \citep*{boumal_voroninski_bandeira}, the additional assumption is the third one. We expect it to be satisfied in most applications: A linear map between two vector spaces is generically injective when the dimension of the second space is at least as large as the dimension of the first one, so when $p(p+1)/2 \leq m$, the map $\phi_V$ in Definition \ref{def:face_regular} should be a priori injective for ``generic'' $V$.
The proof of Theorem \ref{theorem:bvb_improved} is in Section \ref{s:proof_bvb_improved}.

\subsection{Regime where there may be bad critical points\label{ss:main}}

We can now address our main question: How optimal is the result of the previous section? When Problem \eqref{eq:SDP} has a unique global minimizer, with rank $r$ of the order of $\sqrt{2m}$, the result cannot be significantly improved: $p\geq r$ is a necessary condition for Problems \eqref{eq:SDP} and \eqref{eq:P} to have the same minimum. However, as said in the introduction, Problem \eqref{eq:SDP} often admits a solution with rank $r\ll\sqrt{2m}$, and the Burer-Monteiro factorization is numerically observed to work when $p=O(r)$.

Our main theorem however states that, even if we assume $r\ll \sqrt{2m}$, our previous result is essentially not improvable without additional hypotheses on $C$: Under reasonable assumptions on $(\mathcal{A},b)$, if $\frac{p(p+1)}{2} + pr\leq m$, there is a set of cost matrices with non-zero Lebesgue measure for which Problem \eqref{eq:SDP} has a rank $r$ optimal solution, but Problem \eqref{eq:P} has a non-optimal second-order critical point. In particular, for $r=1$, the inequality $\frac{p(p+1)}{2}+p>m$ in Theorem \ref{theorem:bvb_improved} is exactly optimal.

\begin{thm}\label{theorem:main}
  Let $r\in\mathbb{N}^*$ be fixed. Let $p \geq r$ be such that
  \begin{equation*}
    \frac{p(p+1)}{2}+pr \leq m.
  \end{equation*}
  We make the following hypotheses:
  \begin{enumerate}
  \item $\mathcal{C}$ has at least one extreme point with rank $r$, denoted by $X_0$;
  \item $(\mathcal{A},b)$ is $p$-regular;
  \item There exists $V\in\mathcal{M}_p$ such that $\mathcal{M}_p$ is $X_0$-minimally secant at $V$.
  \end{enumerate}
  
  Then there exists a subset $\mathcal{E}_{bad}$ of $\mathbb{S}^{n\times n}$ with non-zero Lebesgue measure such that, for any cost matrix $C\in\mathcal{E}_{bad}$,
  \begin{itemize}
  \item Problem \eqref{eq:SDP} has a unique global minimizer, which has rank $r$.
  \item Problem \eqref{eq:P} has at least one second-order critical point that is not a global minimizer.
  \end{itemize}
\end{thm}
The proof of Theorem \ref{theorem:main} is in Section \ref{s:proof}.

\begin{rem}
  Theorem \ref{theorem:main} stays valid if one replaces ``second-order critical point'' with ``local minimizer''. Indeed, it turns out that the second-order critical points constructed in our proof are non-degenerate, and therefore, local minimizers.
\end{rem}

\begin{rem}\label{rem:non_tight} 
  The inequalities $\frac{p(p+1)}{2}+p>m$ and $\frac{p(p+1)}{2}+pr\leq m$ in Theorems \ref{theorem:bvb_improved} and \ref{theorem:main} are exactly complementary when $r=1$. When $r\geq 2$, there is a small gap between them. Appendix \ref{app:non_tight} shows through an example that there are settings where $\frac{p(p+1)}{2}+p\leq m<\frac{p(p+1)}{2}+pr$ and the conclusions of Theorem \ref{theorem:main} still hold, but we do not know whether it is always the case.
\end{rem}

\subsection{Examples\label{ss:examples}}

\subsubsection{MaxCut\label{sss:maxcut}}
In this subsection, we apply our results to the most famous instance of a problem with the form \eqref{eq:SDP}, the \textit{MaxCut} relaxation:
\begin{align}
  \mbox{minimize } &\scal{C}{X}\nonumber\\
  \mbox{such that }&\mathrm{diag}(X)=1, \label{eq:SDP-MaxCut}\tag{SDP-Maxcut}\\
                   &X\succeq 0.\nonumber
\end{align}

This problem is a relaxation of the ``maximum cut'' problem from graph theory \citep*{delo93, poljak_rendl}, made famous by the work \citep*{goem95}. It also appears in phase retrieval \citep*{maxcut} and $\mathbb{Z}/2\mathbb{Z}$ synchronization \citep*{abbe,bandeira_low_rank} (in which cases its global optimizer is known, both theoretically and numerically, to often have very low rank, typically $1$).

Theorems \ref{theorem:bvb_improved} and \ref{theorem:main} exactly describe when its Burer-Monteiro factorization has no non-optimal second-order critical point for almost any cost matrix, even if we assume that the global minimizer has rank $1$.

\begin{cor}\label{cor:maxcut}
  If $p\in\mathbb{N}$ is such that $
    \frac{p(p+1)}{2}+p > n $
  then, for almost any cost matrix $C$, all second-order critical points of the Burer-Monteiro factorization of Problem \eqref{eq:SDP-MaxCut} are globally optimal.  
    
  On the other hand, for any $p$ such that $
    \frac{p(p+1)}{2} + p \leq n $,
  the set of cost matrices admits a subset with non-zero Lebesgue measure on which
  \begin{itemize}
  \item Problem \eqref{eq:SDP-MaxCut} has a unique global minimizer, which has rank $1$;
  \item Its Burer-Monteiro factorization with rank $p$ has at least one non-optimal second-order critical point.
  \end{itemize}
\end{cor}

This result is proved in Appendix \ref{app:maxcut}.

\subsubsection{Orthogonal-Cut\label{sss:orthogonal_cut}}

We now consider a generalization of \textit{MaxCut}, coined \textit{Orthogonal-Cut} in \citep*{bandeira_kennedy_singer}:
\begin{align}
  \mbox{minimize } &\scal{C}{X}\nonumber\\
  \mbox{such that }&X\in\mathbb{S}^{Sd \times Sd},\nonumber\\
                   & \mathrm{Block}_s(X)=I_d,\forall s=1,\dots,S,\label{eq:SDP-Orthogonal-Cut}\tag{SDP-Orthogonal-Cut}\\
                   &X\succeq 0,\nonumber
\end{align}
where $d,S$ belong to $\mathbb{N}^*$ (with, typically, $d=1,2$ or $3$) and, for any $M\in\R^{Sd\times Sd}$, $s\leq S$, $\mathrm{Block}_s(M)$ is the $s$-th diagonal $d\times d$ block of $M$. Observe that this is exactly Problem \eqref{eq:SDP-MaxCut} when $d=1$.

Problem \eqref{eq:SDP-Orthogonal-Cut} is a natural relaxation of some optimization problems on $O(d)^S$. It notably has applications in molecular imaging \citep*{wang_singer_wen}, sensor network localization \citep*{cucuringu_lipman_singer} and ranking \citep*{cucuringu}. For some theoretical analysis of this semidefinite problem, including conditions under which it admits a low-rank global minimizer, the reader can refer, not only to \citep*{bandeira_kennedy_singer}, but to \citep*{chaudhury_khoo_singer}, \citep*{rosen_carlone_bandeira_leonard} or \citep*{eriksson_olsson_kahl_chin}, as well.

Problem \eqref{eq:SDP-Orthogonal-Cut} is exactly equivalent to
\begin{align*}
  \mbox{minimize } &\scal{C}{X}\\
  \mbox{such that }&X\in\mathbb{S}^{Sd\times Sd},\\
                   &\mathcal{A}(X)=b,\\
                   &X\succeq 0,
\end{align*}
with
$\mathcal{A} : X\in\mathbb{S}^{Sd\times Sd} \to (T_{sup}(\mathrm{Block}_1(X)), \dots, T_{sup}(\mathrm{Block}_S(X))) \in \R^{Sd(d+1)/2}$,
and
$b = (T_{sup}(I_d),\dots,T_{sup}(I_d)) \in \R^{Sd(d+1)/2}$,
where $T_{sup}:\R^{d\times d}\to \R^{d(d+1)/2}$ is the operator that extracts the $\frac{d(d+1)}{2}$ coefficients of the upper triangular part of a matrix.

With these definitions, $(\mathcal{A},b)$ is $p$-regular for any $p\in\mathbb{N}$. In particular, $\mathcal{M}_p$ is a manifold (non-empty if and only if $p\geq d$).
\begin{cor}\label{cor:orthogonal_cut}
  Let us assume that $d=1,2$ or $3$. If
  \begin{equation*}
    \frac{p(p+1)}{2} + p > \frac{Sd(d+1)}{2},
  \end{equation*}
  then, for almost any cost matrix $C$, all second-order critical points of the Burer-Monteiro factorization of Problem \eqref{eq:SDP-Orthogonal-Cut} are globally optimal.
  
  On the other hand, for any $p\geq d$ such that
  \begin{equation*}
    \frac{p(p+1)}{2} + pd \leq \frac{Sd(d+1)}{2},
  \end{equation*}
  there is a set of cost matrices, with non-zero Lebesgue measure, on which
  \begin{itemize}
  \item Problem \eqref{eq:SDP-Orthogonal-Cut} has a unique minimizer, whose rank is $d$;
  \item Its Burer-Monteiro factorization with rank $p$ has at least one non-optimal second-order critical point.
  \end{itemize}
\end{cor}

The proof is in Appendix \ref{app:orthogonal_cut}.

\subsubsection{Optimization over a product of spheres\label{sss:product_of_spheres}}

As a final example, let us consider the problem
\begin{align}
  \mbox{minimize } &\scal{C}{X}\nonumber\\
  \mbox{such that }&X\in\mathbb{S}^{D \times D},\nonumber\\
                   &\sum_{k=d_1+\dots+d_{s-1}+1}^{d_1+\dots+d_s}X_{k,k}=1 ,\forall s=1,\dots,S,\label{eq:SDP-Product}\tag{SDP-Product}\\
                   &X\succeq 0,\nonumber
\end{align}
where $S,d_1,\dots,d_S$ belong to $\mathbb{N}^*$, and $D=d_1+\dots+d_S$. This is the natural semidefinite relaxation of problems that consist in minimizing a degree $2$ polynomial function on the product of spheres $S^{d_1}\times \dots \times S^{d_S}$.

Problem \eqref{eq:SDP-Product} encompasses several important particular cases: when $d_1=\dots=d_S=1$, we recover Problem \eqref{eq:SDP-MaxCut}. When $d_1=\dots=d_S=2$, it is equivalent to a complex version of \eqref{eq:SDP-MaxCut} (for matrices $C$ of a particular form). When $S=2$ and $d_2=1$, it is the relaxation of a standard trust-region subproblem \citep*[Subsection 5.2]{boumal_voroninski_bandeira}. For general values of $d_1,\dots,d_S$, it is a simplification of the relaxation of optimization problems over an intersection of ellipsoids, which appear in trust-region algorithms for constrained problems \citep*{celis}.

\begin{cor}\label{cor:product_of_spheres}
  If $p\in\mathbb{N}$ is such that
  \begin{equation*}
    \frac{p(p+1)}{2} + p > S,
  \end{equation*}
  then, for almost any cost matrix $C$, all second-order critical points of the Burer-Monteiro factorization of Problem \eqref{eq:SDP-Product} are globally optimal.
  
  On the other hand, for any $p\in\mathbb{N}^*$ such that
  \begin{equation*}
    \frac{p(p+1)}{2} + p \leq S,
  \end{equation*}
  the set of cost matrices admits a subset with non-zero Lebesgue measure on which
  \begin{itemize}
  \item Problem \eqref{eq:SDP-Product} has a unique global optimum, which has rank $1$;
  \item Its Burer-Monteiro factorization with rank $p$ has at least one non-optimal second-order critical point.
  \end{itemize}  
\end{cor}

The proof is in Appendix \ref{app:product_of_spheres}.

\section{Proof of Theorem \ref{theorem:bvb_improved}\label{s:proof_bvb_improved}}

Let $m$ be such that
\begin{equation}\label{eq:bvb_improved_value_of_m}
  m < \frac{p(p+1)}{2}+p.
\end{equation}

From \citep*[Theorem 2.1]{pataki}, Problem \eqref{eq:SDP} has a minimizer with rank at most $p$. Consequently, Problems \eqref{eq:SDP} and \eqref{eq:P} have the same minimum, and if $V$ is a global minimizer of Problem \eqref{eq:P}, $X=VV^T$ is a minimizer of Problem \eqref{eq:SDP}. It therefore suffices to show that, for almost all cost matrices, Problem \eqref{eq:P} has no second-order critical point which is not a global minimizer.

  \subsection{Overview of the proof}

  Our proof starts in a similar way as the one in \citep*{boumal_voroninski_bandeira}. Namely, we use the first and second-order properties of critical points to parametrize the set of ``bad'' cost matrices: we define an appropriate manifold $\mathcal{M}_{param}$ and a smooth map $\phi:\mathcal{M}_{param}\to\mathbb{S}^{n\times n}$ such that the set of bad cost matrices is included in $\phi(\mathcal{M}_{param})$.

  Then the proofs differ. The authors of \citep*{boumal_voroninski_bandeira} show that, when $\frac{p(p+1)}{2} > m$, the dimension of their manifold $\mathcal{M}_{param}$ is strictly smaller than $\dim(\mathbb{S}^{n\times n})$, hence $\phi(\mathcal{M}_{param})$ has zero Lebesgue measure in $\mathbb{S}^{n\times n}$. On our side, we use additional properties of critical points to show that the set of bad cost matrices is actually included in the critical values of $\phi$, and not only in the range of $\phi$. This set has zero Lebesgue measure in $\mathbb{S}^{n\times n}$, from Sard's theorem.

  \subsection{Details}
  
  We define $\mathcal{M}_p^{full}=\{V\in\mathcal{M}_p,\mathrm{rank}(V)=p\}$. It is an open subset of $\mathcal{M}_p$, and therefore also a $(np-m)$-dimensional Riemannian manifold. We also define
  \begin{equation*}
    \mathcal{E} = \{(V,C_2)\in\mathcal{M}_p^{full}\times \mathbb{S}^{n\times n}\mbox{ such that }C_2V=0_{n,p}\}.
  \end{equation*}
  This set is a manifold, as stated in the following proposition, whose proof is in Appendix \ref{app:mathcalE}.
  \begin{prop}\label{prop:mathcalE}
    The set $\mathcal{E}$ is a manifold, with dimension $np-m+\frac{(n-p)(n-p+1)}{2}$. Additionally, for any $(V,C_2)\in\mathcal{E}$,
    \begin{equation}\label{eq:tangent_mathcalE}
      T_{(V,C_2)}\mathcal{E}
       = \{(\dot V,\dot C_2)\in T_V\mathcal{M}_p \times \mathbb{S}^{n\times n}\mbox{ such that } \dot C_2 V + C_2\dot V = 0_{n,p}\}.
    \end{equation}
  \end{prop}
  
  We define
  \begin{equation*}
    \begin{array}{cccc}
      \phi:& \mathcal{E} \times \R^m & \to & \mathbb{S}^{n\times n} \\
      & ((V,C_2),\mu) & \to & C_2 + \mathcal{A}^*(\mu).
    \end{array}
  \end{equation*}

  The following lemma, whose proof is in Subsection \ref{ss:critical_value}, says that any cost matrix for which a non-optimal second-order critical point exists is a critical value of $\phi$.
  \begin{lem}\label{lem:critical_value}
    For any cost matrix $C\in\mathbb{S}^{n\times n}$, if Problem \eqref{eq:P} has a non-optimal second-order critical point, then there exist $(V,C_2)\in\mathcal{E},\mu\in\R^m$ such that
    \begin{equation*}
      C = \phi((V,C_2),\mu),
    \end{equation*}
    and the mapping $d\phi((V,C_2),\mu):T_{(V,C_2)}\mathcal{E} \times \R^m \to \mathbb{S}^{n\times n}$ is not surjective.
  \end{lem}
  From Sard's theorem, we can therefore conclude that the set of such cost matrices has zero measure in $\mathbb{S}^{n\times n}$.

  \subsection{Proof of Lemma \ref{lem:critical_value}\label{ss:critical_value}}

  Let $C\in\mathbb{S}^{n\times n}$ be a cost matrix for which a non-optimal second-order critical point exists. Let $V\in\mathcal{M}_p$ be such a critical point. From \citep*[Theorem 1.6]{boumal_voroninski_bandeira}, $\mathrm{rank}(V)=p$, so $V$ is in $\mathcal{M}_p^{full}$.

  As $V$ is first-order critical, there exist, from Proposition \ref{prop:first_order}, $\mu\in\R^m,C_2\in\mathbb{S}^{n\times n}$ such that $C_2V=0$ and
  %\citep*[Eqs (8) and (11)]{boumal_voroninski_bandeira}, $\mu\in\R^m$ such that $(C-\mathcal{A}^*(\mu))V=0$. Setting $C_2=C-\mathcal{A}^*(\mu)$, we see that
  \begin{equation*}
    C = C_2 + \mathcal{A}^*(\mu) = \phi((V,C_2),\mu).
  \end{equation*}

  Let us now show that $d\phi((V,C_2),\mu)$ is not surjective. From \citep*[Theorem 1.6]{boumal_voroninski_bandeira}, the dimension of the face of $\mathcal{C}$ containing $VV^T$ is at least  
  \begin{equation*}
    \frac{p(p+1)}{2}-m+p \overset{\textrm{Eq. \eqref{eq:bvb_improved_value_of_m}}}> 0.
  \end{equation*}
  Let $X_{Face}\ne VV^T$ be an element of this face. Using the geometrical properties of $\{X\in\mathbb{S}^{n\times n},X\succeq 0\}$, one can establish the following proposition, whose proof is in Appendix \ref{app:face_Splus}.
  \begin{prop}\label{prop:face_Splus}
    There exists $T\in\mathbb{S}^{p\times p}$ such $X_{Face}=VTV^T$.
  \end{prop}
  Let $T$ be as in the proposition.
  For any $(\dot V,\dot C_2)\in T_{(V,C_2)}\mathcal{E},\dot\mu\in\R^m$,
  \begin{align*}
    \langle d\phi((V,C_2),\mu)& \cdot((\dot V,\dot C_2),\dot\mu),
      X_{Face}-VV^T\rangle \\
    & = \scal{\dot C_2 + \mathcal{A}^*(\dot\mu)}{X_{Face}-VV^T} \\
    & = \scal{\dot C_2}{V(T-I_p)V^T} + \scal{\dot\mu}{\mathcal{A}(X_{Face}-VV^T)} \\
    & \overset{(a)}{=} \scal{\dot C_2}{V(T-I_p)V^T} \\
    & = \scal{\dot C_2 V}{V(T-I_p)} \\
    & \overset{(b)}{=} - \scal{C_2\dot V}{V(T-I_p)} \\
    & = -\scal{\dot V}{C_2 V (T-I_p)} \\
    & \overset{(c)}{=} 0.
  \end{align*}
  Equality $(a)$ is true because $X_{face}$ belongs to $\mathcal{C}$, so $\mathcal{A}(X_{Face})=b=\mathcal{A}(VV^T)$. Equality $(b)$ is true because of Equation \eqref{eq:tangent_mathcalE}, and equality $(c)$ because $C_2 V=0_{n,p}$ from the definition of $\mathcal{E}$.

  This shows that the range of $d\phi((V,C_2),\mu)$ in $\mathbb{S}^{n\times n}$ is included in $(X_{Face}-VV^T)^\perp$, so that $d\phi((V,C_2),\mu)$ cannot be surjective.
    
\section{Proof of Theorem \ref{theorem:main}\label{s:proof}}

This section is devoted to the proof of the main theorem. The first two subsections, \ref{ss:main_first_part} and \ref{ss:main_second_part}, each present the outline of one half of the proof, with the technical details hidden into lemmas. The remaining subsections contain the proofs of these lemmas.

\subsection{First part\label{ss:main_first_part}}

In the first part, we assume (proving this assumption is done in the second part) that there exists one cost matrix, $C$, for which Problem \eqref{eq:SDP} has a unique global minimizer, with rank $r$, but Problem \eqref{eq:P} has a spurious second-order critical point. We show that, for all matrices close enough to $C$, these properties still hold, hence they hold on a whole set, with non-zero Lebesgue measure.

As stated, this assertion may not be quite true ($C$ might be an isolated ``bad'' cost matrix). However, it becomes true if we assume $C$ to satisfy some additional non-degeneracy properties.
Consequently, in this part of the proof, we admit the following lemma.
\begin{lem}\label{lem:one_bad_C}
  There exists a cost matrix $C\in\mathbb{S}^{n\times n}$ such that
  \begin{itemize}
  \item Problem \eqref{eq:SDP} has a unique global minimizer, whose rank is $r$.
  \item Strict complementary slackness holds.
  \item Problem \eqref{eq:P} has a second-order critical point, which is not a global minimizer.
  \item This second-order critical point is non-degenerate.
  \end{itemize}
\end{lem}

We respectively denote $X_0$ and $V$ the global minimizer and second-order critical point of Lemma \ref{lem:one_bad_C}.
The two properties we must show are stated in the following lemmas.
\begin{lem} \label{lem:global_min_SDP}
  For any matrix $C'$ close enough to $C$, Problem \eqref{eq:SDP} has a unique global minimizer, and this minimizer has rank $r$.
\end{lem}

\begin{lem}\label{lem:second_order_stability}
  For any matrix $C'$ close enough to $C$, Problem \eqref{eq:P} has a second-order critical point which is not a global minimizer.
\end{lem}

To prove Lemma \ref{lem:global_min_SDP}, we use general convexity and continuity arguments to show that at least one minimizer exists, and that it goes to $X_0$ when $C'$ goes to $C$. In particular, it has rank at least $r$ when $C'$ is close enough to $C$. With another continuity argument, we show that, because strict complementary slackness holds for $C$, it also holds for any $C'$ close enough to $C$. Therefore, the minimizer is unique (from Proposition \ref{prop:strict_complementary_slackness}), and strict complementary slackness also allows us to prove that it has rank exactly $r$. A detailed proof is in Subsection \ref{ss:global_min_SDP}.

For Lemma \ref{lem:second_order_stability}, it actually suffices to show that, for $C'$ close to $C$, Problem \eqref{eq:P} has a second-order critical point close to $V$. Indeed, no matrix close enough to $V$ can be a global minimizer. If the Hessian at $V$ was positive definite, this would follow from general geometric arguments. But because of the invariance of the problem to multiplication by elements of $O(p)$, the Hessian is not positive definite. We must therefore consider the quotient manifold $\mathcal{M}_p/O(p)$ and a quotiented version of Problem \eqref{eq:P}. For this version, as $V$ is non-degenerate, the Hessian is positive definite, so the general arguments apply.
The details are in Subsection \ref{ss:second_order_stability}.

\subsection{Second part: proof of lemma \ref{lem:one_bad_C}\label{ss:main_second_part}}

We recall that we want to construct a cost matrix $C$ such that
\begin{enumerate}
\item Problem \eqref{eq:SDP} has a unique global minimizer, with rank $r$.
\item Strict complementary slackness holds.
\item Problem \eqref{eq:P} has a second-order critical point, which is not a global minimizer.
\item This second-order critical point is non-degenerate.
\end{enumerate}
It turns out that, for any rank $r$ matrix $X_0\in\mathcal{C}$ and any $V\in\mathcal{M}_p$, provided that $X_0$ is extremal in $\mathcal{C}$ and $\mathcal{M}_p$ is $X_0$-minimally secant at $V$, it is possible to construct a matrix $C$ as desired, and such that in addition, the unique global minimizer is precisely $X_0$, and the spurious critical point is precisely $V$.

Let us fix $X_0,V$ as described (we have made the hypothesis they existed) and explain how to construct $C$.
First, the results in Subsection \ref{ss:basic} allow us to rephrase the desired conditions in more analytical terms: they are equivalent to the existence of $g_1,g_2\in\mathbb{R}^m,C_1,C_2\in\mathbb{S}^{n\times n}$ such that
\begin{subequations}
\begin{align}
  C & = \mathcal{A}^*(g_1)+C_1 ; \label{eq:12_cond1}\\
  C_1 & \succeq 0 ; \label{eq:12_cond2}\\
  C_1 X_0 & = 0 ; \label{eq:12_cond3}\\
  \mathrm{rank}(C_1) & = n-\mathrm{rank}(X_0)=n-r ;\label{eq:12_cond4} \\
  \mathcal{A}^*(g_1)+C_1 & = \mathcal{A}^*(g_2)+C_2 ; \label{eq:34_cond1}\\
  C_2 V & = 0 ; \label{eq:34_cond2}\\
  \forall \dot V\in T_V\mathcal{M}_p,\quad
  \scal{C_2}{\dot V \dot V^T} & \geq 0,\mbox{ with equality}\nonumber\\
    & \quad \mbox{ iff }\dot V = VA, A\in\Anti(p).\label{eq:34_cond3}
\end{align}
\end{subequations}

The construction now proceeds as follows:
\begin{enumerate}
\item We set $g_2=0$.\label{it:construction1}
\item We construct $g_1,C_1,C_2$ such that Properties \eqref{eq:12_cond2}, \eqref{eq:12_cond3}, \eqref{eq:12_cond4}, \eqref{eq:34_cond1} and \eqref{eq:34_cond2} hold.\label{it:construction2}
\item From $g_1,C_1,C_2$, we construct $g_1^{(mod)},C_1^{(mod)},C_2^{(mod)}$ which satisfy Property \eqref{eq:34_cond3} in addition to the previous ones.\label{it:construction3}
\item We set $C=\mathcal{A}^*(g_1^{(mod)})+C_1^{(mod)}$; it satisfies all the required properties.
  \label{it:construction4}
\end{enumerate}
Points \ref{it:construction1} and \ref{it:construction4} are straightforward. For Points \ref{it:construction2} and \ref{it:construction3}, see Subsections \ref{ss:point2} and \ref{ss:point3}.

\subsection{Proof of Lemma \ref{lem:global_min_SDP}\label{ss:global_min_SDP}}
  
% \begin{lem*}[Lemma \ref{lem:global_min_SDP}]
%   For some cost matrix $C$, we assume that Problem \eqref{eq:SDP} has $X_0$ as unique global minimizer. We also assume that $C$ can be written as $
%     C = C_1 + \mathcal{A}^*(g_1),$
%   with $g_1\in\R^m$ and $C_1\in\mathbb{S}^{n\times n}$ such that $C_1\succeq 0,\mathrm{rank}(C_1)=n-r$ and $C_1X_0=0$.
  
%   Then, for any matrix $C'$ close enough to $C$, Problem \eqref{eq:SDP} (with cost matrix $C'$) also has a unique global minimizer, and this minimizer has rank $r$.
% \end{lem*}

% \begin{proof}
To establish the lemma, it suffices to show that, for any sequence $(C'_k)_{k\in\mathbb{N}}$ of cost matrices converging to $C$, Problem \eqref{eq:SDP} with cost matrix $C'_k$ has a unique minimizer, and this minimizer has rank $r$, as soon as $k$ is large enough. Let $(C'_k)_{k\in\mathbb{N}}$ be such a sequence. 

The following proposition shows that, for $k$ large enough, at least one minimizer exists, and it is arbitrarily close to $X_0$.
Its proof is in Appendix \ref{app:argmin_close_to_X0}.
  \begin{prop}\label{prop:argmin_close_to_X0}
    Let $\epsilon>0$ be fixed. For $k$ large enough,
    \begin{itemize}
    \item Problem \eqref{eq:SDP} (with cost matrix $C'_k$) admits at least one minimizer;
    \item all minimizers of Problem \eqref{eq:SDP} belong to the ball $B(X_0,\epsilon)$.
    \end{itemize}
  \end{prop}
  For any $k$ large enough, let $X'_k$ be a minimizer corresponding to the cost matrix $C'_k$. If there are several of them, we choose $X'_k$ as an extremal point of the set of minimizers (such a point exists because the set is bounded, from Proposition \ref{prop:argmin_close_to_X0}, convex and closed); it is then also an extremal point of the feasible set $\mathcal{C}$.
Let us show that, for $k$ large enough,
  \begin{equation}\label{eq:global_min_sequence}
    \mathrm{rank}(X'_k)=r\mbox{ and }X'_k\mbox{ is the unique minimizer of Problem \eqref{eq:SDP}.}
  \end{equation}

  Let $g_1,C_1$ be defined as in Proposition \ref{prop:KKT_SDP}: $C=\mathcal{A}^*(g_1)+C_1,C_1\succeq 0$ and $C_1X_0=0$. Similarly, let, for any $k$, $h_k\in\R^m,D_k\in\mathbb{S}^{n\times n}$ be such that $C'_k=\mathcal{A}^*(h_k)+D_k,D_k\succeq 0$ and $D_k X'_k = 0$.
  
  The following lemma states that $D_k\overset{k\to+\infty}{\to}C_1$ and $h_k\overset{k\to+\infty}{\to}{g_1}$. Its proof is in Subsection \ref{app:convergence_DC_hg} and relies on the $p$-regularity of $(\mathcal{A},b)$.
  \begin{lem}\label{lem:convergence_DC_hg}
    When $k$ goes to infinity, $D_k \to C_1$ and $h_k\to g_1.$
  \end{lem}
  
  For any $k$, because $D_kX'_k=0$,
  \begin{equation}\label{eq:sum_ranks_D_X}
    \mathrm{rank}(D_k) + \mathrm{rank}(X'_k) \leq n.
  \end{equation}
  From Proposition \ref{prop:argmin_close_to_X0}, $(X'_k)_{k\in\mathbb{N}}$ converges to $X_0$, and from Lemma \ref{lem:convergence_DC_hg}, $(D_k)_{k\in\mathbb{N}}$ converges to $C_1$. In particular, for $k$ large enough,
  \begin{gather*}
    \mathrm{rank}(X'_k) \geq \mathrm{rank}(X_0) = r \\
    \mbox{and}\quad \mathrm{rank}(D_k)\geq \mathrm{rank}(C_1) = n-r.
  \end{gather*}
  Combined with Equation \eqref{eq:sum_ranks_D_X}, this proves that, for $k$ large enough, $
    \mathrm{rank}(X'_k)=r\quad\mbox{and}\quad \mathrm{rank}(D_k)=n-r. $
  This establishes the first part of Property \eqref{eq:global_min_sequence}. The second part is a direct consequence of Proposition \ref{prop:strict_complementary_slackness}.
%\end{proof}
  
  \subsection{Proof of Lemma \ref{lem:second_order_stability}\label{ss:second_order_stability}}

  We recall that $X_0$ is a global minimizer of Problem \eqref{eq:SDP}, but $VV^T$ is not: $\scal{C}{X_0}<\scal{C}{VV^T}$. By continuity, there is actually a neighborhood $\mathcal{V}$ of $V$ in $\mathcal{M}_p$ such that
  \begin{equation}\label{eq:value_at_X0_Vprime}
    \forall V'\in \mathcal{V},\quad \scal{C}{X_0} < \scal{C}{V'V'^T}.
  \end{equation}
  By continuity again, Equation \eqref{eq:value_at_X0_Vprime} stays true if we replace $C$ by any close enough matrix $C'$. Therefore, for $C'$ close enough to $C$, no matrix of the form $V'V'^T$ with $V'\in\mathcal{V}$ can be a global minimizer of Problem \eqref{eq:SDP}, hence no matrix of $\mathcal{V}$ can be a global minimizer of Problem \eqref{eq:P}. From this remark, if we show that, for any $C'$ close enough to $C$, Problem \eqref{eq:P} has a second-order critical point in $\mathcal{V}$, we have proved the lemma. Let us do that.
  
  For any cost matrix $C'\in\mathbb{S}^{n\times n}$, we denote
  $f_{C'}:W\in\mathcal{M}_p \to \scal{C'}{WW^T} \in\R$ the cost function of Problem \eqref{eq:P}.
  If $\mathrm{Hess}f_C(V)$ was positive definite, we could apply the following general proposition (proved in Appendix \ref{app:second_order_stability_general}).
  \begin{prop}\label{prop:second_order_stability_general}
    Let $\mathcal{M}$ be a Riemannian manifold, $E$ a finite-dimensional vector space, and $f:E\times \mathcal{M}\to\R$ a smooth map. Let $c\in E,v\in \mathcal{M}$ be fixed. We assume that $f(c,.)$ has a second-order critical point at $v$, and that $\mathrm{Hess}(f(c,.))(v)\succ 0$.

    Then, for any neighborhood $\mathcal{V}$ of $v$ in $\mathcal{M}$, the map $f(c',.)$ has a second-order critical point in $\mathcal{V}$ for any $c'\in E$ close enough to $c$.
  \end{prop}
  However, because $f_C$ is invariant to right multiplication by elements of $O(p)$, the Hessian is degenerate. Therefore, before applying the proposition, we must explicitly factorize this invariance by introducing the corresponding quotient manifold. We refer to \citep*[Section 3.4]{absil} for basic results on quotient manifolds.
Specifically, let $\mathcal{M}_p^{full}$ be the open subset of $\mathcal{M}_p$ that contains its rank $p$ elements\footnote{If $\mathcal{M}_p$ contains rank-deficient matrices, $\mathcal{M}_p/O(p)$ is not a manifold, because $\{(V,VX),V\in\mathcal{M}_p,X\in O(p)\}$ is not a submanifold of $\mathcal{M}_p^2$. We must therefore remove rank-deficient elements from $\mathcal{M}_p$.}. The quotient $\mathcal{M}_p^{full}/O(p)$ is a manifold with dimension
  \begin{equation*}
    \dim(\mathcal{M}_p) - \dim(O(p)) = \dim(\mathcal{M}_p) - \frac{p(p-1)}{2}.
  \end{equation*}
  Since $O(p)$ acts by isometries on $\mathcal{M}_p^{full}$, $\mathcal{M}_p^{full}/O(p)$ inherits from the Riemannian structure of $\mathcal{M}_p^{full}$.
  We denote $Q:\mathcal{M}_p^{full} \to \mathcal{M}_p^{full}/O(p)$ the canonical projection. It is a smooth map, with surjective differential everywhere.

  For any $C'$, since $f_{C'}$ is invariant to the action of $O(p)$, we can define its quotient, that is the (also smooth) map $f_{C',O(p)}:\mathcal{M}_p^{full}/O(p)\to\mathbb{R}$ such that
  % voir prop 3.4.5 d'[Absil et al] si une référence est demandée
  \begin{equation*}
    f_{C',O(p)}\circ Q = f_{C'}.
  \end{equation*}

  The following proposition (proved in Appendix \ref{app:transfer_critical_points}) shows that there is a correspondance between the critical points of $f_{C',O(p)}$ and $f_{C'}$.
  \begin{prop}\label{prop:transfer_critical_points}
    Let $\mathcal{N}_1,\mathcal{N}_2$ be two Riemannian manifolds, and $f:\mathcal{N}_2\to\R$ a smooth function. Let $\phi:\mathcal{N}_1\to\mathcal{N}_2$ be a smooth map with surjective differential at any point of $\mathcal{N}_1$.

    Then, for any $v\in\mathcal{N}_1$, $v$ is a second-order critical point of $f\circ \phi$ if and only if $\phi(v)$ is a second-order critical point of $f$. Additionally,
    \begin{equation*}
      \mathrm{rank}(\mathrm{Hess}(f\circ \phi)(v)) =
      \mathrm{rank}(\mathrm{Hess}f(\phi(v))).
    \end{equation*}
  \end{prop}

  This proposition, applied to $\mathcal{N}_1=\mathcal{M}_p^{full}, \mathcal{N}_2=\mathcal{M}_p^{full}/O(p)$, $f=f_{C',O(p)}$ and $\phi=Q$, shows that, because $V$ is a second-order critical point of $f_{C}=f_{C,O(p)}\circ Q$, $Q(V)$ is a second-order critical point of $f_{C,O(p)}$ and
  \begin{align*}
    \mathrm{rank}(\mathrm{Hess}f_{C,O(p)}(Q(V)))
    & = \mathrm{rank}(\mathrm{Hess}f_C(V)) \\
    & \overset{(Rem. \ref{rem:non_degenerate_rank})}= \dim(\mathcal{M}_p)-\frac{p(p-1)}{2} \\
    & = \dim(\mathcal{M}_p^{full}/O(p)).
  \end{align*}

  In other words, $\mathrm{Hess}f_{C,O(p)}(Q(V))$ is positive definite.

  We apply Proposition \ref{prop:second_order_stability_general} to $E=\mathbb{S}^{n\times n},\mathcal{M}=\mathcal{M}_p^{full}/O(p)$ and $f:(C',W)\in\mathbb{S}^{n\times n}\times\mathcal{M}_p^{full}/O(p) \to f_{C',O(p)}(W)$: for any neighborhood $\mathcal{V}_{O(p)}$ of $Q(V)$, $f_{C',O(p)}$ has a second-order critical point in $\mathcal{V}_{O(p)}$ if $C'$ is close enough to $C$.

  We use this property with $\mathcal{V}_{O(p)} = Q(\mathcal{V})$. For any $C'$ close enough to $C$, $f_{C',O(p)}$ has a second-order critical point of the form $Q(W)$, with $W\in \mathcal{V}$. Then, from Proposition \ref{prop:transfer_critical_points}, $f_{C'}$ has a second-order critical point in $\mathcal{V}$.

\subsection{Construction of $C$: Point \ref{it:construction2}\label{ss:point2}}    

We must show the existence of $g_1,C_1,C_2$ such that
\begin{subequations}
\begin{gather}
  C_1 \succeq 0,\quad C_1X_0=0,\quad \mathrm{rank}(C_1)=n-r,\label{eq:cond12}\\
  \mathcal{A}^*(g_1)+C_1 = C_2,\quad C_2 V=0.\label{eq:cond34}
\end{gather}
\end{subequations}

We simplify the problem with the following proposition, proved in Appendix \ref{app:up_to_replacing}.
\begin{prop}\label{prop:up_to_replacing}
  Without loss of generality, we can assume that
  \begin{equation*}
    X_0 = \left(\begin{smallmatrix}
        I_r&0_{r,n-r}\\
        0_{n-r,r}&0_{n-r,n-r}
      \end{smallmatrix}\right)
    \quad\mbox{and}\quad
    V = \left(\begin{smallmatrix}0_{r,p} \\ I_p \\ 0_{n-p-r,p} \end{smallmatrix}\right).
  \end{equation*}
\end{prop}
With this assumption, the three conditions in Equation \eqref{eq:cond12} are true if and only if
\begin{equation*}
  C_1 = \left(\begin{smallmatrix}0_{r,r}&0_{r,n-r}\\ 0_{n-r,r}&D_1 \end{smallmatrix}\right)
\end{equation*}
for some $D_1\in\mathbb{S}^{(n-r)\times(n-r)}$ such that $D_1\succ 0$. And $C_2V=0$ if and only if
\begin{equation*}
  C_2 = \left(\begin{smallmatrix}F_1&0_{r,p}&F_2\\
      0_{p,r}&0_{p,p}&0_{p,n-r-p}\\
      F_2^T&0_{n-r-p,p}&F_3 \end{smallmatrix}\right),
\end{equation*}
for some $F_1,F_2,F_3$. Therefore, to ensure Conditions \eqref{eq:cond12} and \eqref{eq:cond34}, we must only show the existence of $g_1,D_1,F_1,F_2,F_3$ such that $D_1\succ 0$ and
\begin{equation}\label{eq:cond34_1_equiv}
  \mathcal{A}^*(g_1)
  = \left(\begin{smallmatrix}F_1&0_{r,p}&F_2\\
      0_{p,r}&0_{p,p}&0_{p,n-r-p}\\
      F_2^T&0_{n-r-p,p}&F_3 \end{smallmatrix}\right)
  - \left(\begin{smallmatrix}0_{r,r}&0_{r,n-r}\\ 0_{n-r,r}&D_1 \end{smallmatrix}\right).
\end{equation}

We observe that, if these exist, $\mathcal{A}^*(g_1)$ must be of the form
\begin{equation}\label{eq:cond_on_g1}
  \mathcal{A}^*(g_1)
  = \left(\begin{smallmatrix}G_1&0_{r,p}&G_2
      \\ 0_{p,r}&G_3&G_4 \\ G_2^T&G_4^T&G_5
    \end{smallmatrix}\right),
\end{equation}
with $G_3\prec 0$ (since it is a minor of $-D_1$).
But conversely, if there exists $g_1$ for which Equation \eqref{eq:cond_on_g1} is true, we can set
\begin{equation*}
  F_1=G_1,\quad F_2=G_2,\quad F_3 = G_5 + \lambda I_{n-r-p},
  \quad D_1=\left(\begin{smallmatrix}
      -G_3&-G_4 \\
      -G_4^T&\lambda I_{n-r-p}
    \end{smallmatrix}
  \right),
\end{equation*}
for some $\lambda >0$ large enough, and Equation \eqref{eq:cond34_1_equiv} holds. (We observe that $D_1\succ 0$ for $\lambda$ large enough: all its principal minors are of the form
\begin{equation*}
  \det \left(\begin{smallmatrix} -G_3^{(sub)} & -G_4^{(sub)} \\
      -G_4^{(sub)T}& \lambda I_s
      \end{smallmatrix}\right) = \lambda^s \det(-G_3^{(sub)}) + O(\lambda^{s-1}),
  \end{equation*}
  with $-G_3^{(sub)}$ a principal submatrix of $-G_3$, whose determinant is positive because $-G_3\succ 0$. Therefore, all principal minors of $D_1$ are positive if $\lambda$ is large enough.)

To conclude, we must only prove the existence of $g_1$ for which Equation \eqref{eq:cond_on_g1} is true. This is a consequence of the following lemma, whose proof is in Appendix \ref{app:construction_g1} (and relies on the minimally secant property).
\begin{lem}\label{lem:construction_g1}
  For any $R_1\in\mathbb{R}^{r\times p},R_2\in\mathbb{S}^{p\times p}$,
  there exist $g_1\in\mathbb{R}^m,G_1,G_2,G_4,G_5$ such that,
  \begin{equation*}
    \mathcal{A}^*(g_1)
    = \left(\begin{smallmatrix}G_1&R_1&G_2
        \\ R_1^T&R_2&G_4 \\ G_2^T&G_4^T&G_5
      \end{smallmatrix}\right).
  \end{equation*}
\end{lem}
  
\subsection{Construction of $C$: Point \ref{it:construction3}\label{ss:point3}}

In this subsection, we consider $g_1,C_1,C_2$ satisfying Properties \eqref{eq:12_cond2} to \eqref{eq:34_cond2} and construct $g_1^{(mod)},C_1^{(mod)},C_2^{(mod)}$ which also satisfy these properties, and, in addition, Property \eqref{eq:34_cond3}:
\begin{equation}\label{eq:34_cond3_bis}
  \forall \dot V\in T_V\mathcal{M}_p,\quad
  \scal{C_2^{(mod)}}{\dot V\dot V^T}\geq 0,
\end{equation}
with equality if and only if $\dot V=VA$ for some $A\in\Anti(p)$.

Using Proposition \ref{prop:up_to_replacing} as in the previous subsection, we assume
\begin{equation}\label{eq:up_to_replacing}
  X_0 = \left(\begin{smallmatrix}
      I_r&0_{r,n-r}\\
      0_{n-r,r}&0_{n-r,n-r}
    \end{smallmatrix}\right)
  \quad\mbox{and}\quad
  V = \left(\begin{smallmatrix}0_{r,p} \\ I_p \\ 0_{n-p-r,p} \end{smallmatrix}\right).
\end{equation}

We set
\begin{gather*}
  g_1^{(mod)}=g_1,\quad
  C_1^{(mod)}=C_1 + t \left(\begin{smallmatrix}
      0_{r+p,r+p}&0_{r+p,n-r-p}\\
      0_{n-r-p,r+p}&I_{n-r-p}
    \end{smallmatrix}\right), \\
  C_2^{(mod)}=C_2 + t \left(\begin{smallmatrix}
      0_{r+p,r+p}&0_{r+p,n-r-p}\\
      0_{n-r-p,r+p}&I_{n-r-p}
    \end{smallmatrix}\right),
\end{gather*}
for some $t\geq 0$ large. From the following proposition (proved in Appendix \ref{app:properties_still_hold}) these definitions satisfy Properties \eqref{eq:12_cond2} to \eqref{eq:34_cond2}.
\begin{prop}\label{prop:properties_still_hold}
Whatever the value of $t\geq 0$, $g_1^{(mod)},C_1^{(mod)}$ and $C_2^{(mod)}$ satisfy Properties \eqref{eq:12_cond2} to \eqref{eq:34_cond2}.
\end{prop}

When $t$ is large enough, it turns out that they also satisfy Equation \eqref{eq:34_cond3_bis}. This is proved in two steps, each embedded in a proposition (proofs are in Appendices \ref{app:34_cond3_orth1} and \ref{app:34_cond3_orth2}): first, we observe that, to prove Equation \eqref{eq:34_cond3_bis}, one can look only at matrices $\dot V$ in some subspace of $T_V\mathcal{M}_p$. Then, with a compactness argument, we show that, for matrices in this subspace, Equation \eqref{eq:34_cond3_bis} is true.
\begin{prop}\label{prop:34_cond3_orth1}
  Let $\mathcal{E}_{\perp}$ be the orthogonal in $T_V\mathcal{M}_p$ of $\{VA,A\in\Anti(p)\}$. Equation \eqref{eq:34_cond3_bis} is true if and only if
  \begin{equation}\label{eq:34_cond3_bis_orth}
  \forall \dot V \in\mathcal{E}_{\perp}-\{0\},\quad
  \scal{C_2^{(mod)}}{\dot V \dot V^T} > 0.
\end{equation}
\end{prop}
\begin{prop}\label{prop:34_cond3_orth2}
 For $t$ large enough, Equation \eqref{eq:34_cond3_bis_orth} is true.
\end{prop}

\section*{Acknowledgements}

We thank the reviewers for their careful reading and suggestions, which allowed us to significantly simplify our results and proofs.

\appendix

\section{Proof of basic properties}

\subsection{Proof of Proposition \ref{prop:KKT_SDP}\label{app:KKT_SDP}}

The dual of Problem \eqref{eq:SDP} is
\begin{align*}
  \mbox{maximize }&\scal{g_1}{b}\tag{SDP-dual}\\
  \mbox{such that }&C = \mathcal{A}^*(g_1)+C_1,\\
                  &C_1\succeq 0.
\end{align*}
If $C_1,g_1$ are as in the statement, they are dual feasible. Because of the complementary slackness condition $C_1X_0=0$, $X_0$ and $(C_1,g_1)$ are primal-dual optimal. In particular, $X_0$ is a solution of Problem \eqref{eq:SDP}.

Conversely, let us assume $X_0$ is a solution of Problem \eqref{eq:SDP}. We temporarily admit that Slater's condition holds (that is, $\mathcal{C}$ contains a positive definite matrix). Then strong duality holds \citep*[page 114]{wolkowicz} and the dual problem has at least one solution $(C_1,g_1)$. This pair satisfies $C=\mathcal{A}^*(g_1)+C_1$ and $C_1\succeq 0$, because it is dual feasible. Strong duality means that
\begin{equation*}
\scal{g_1}{b} = \scal{C}{X_0},
\end{equation*}
which is equivalent to $\scal{C_1}{X_0}=0$ and in turn implies $C_1X_0=0$ because $C_1,X_0\succeq 0$.

To establish Slater's condition, we assume by contradiction that it does not hold:
    \begin{equation*}
      \{X\in\mathbb{S}^{n\times n},\mathcal{A}(X)=b\} \cap \{X\in\mathbb{S}^{n\times n},X\succ 0\} = \emptyset.
    \end{equation*}
    From a hyperplane separation theorem, there exists a non-zero $M\in\mathbb{S}^{n\times n}$, and $\mu\in\R$ such that
    \begin{subequations}
      \begin{gather}
        \forall X \in\{X\in\mathbb{S}^{n\times n},X\succ 0\}, \quad \scal{M}{X} > \mu
        \label{eq:separation_2}\\
        \mbox{and}\quad \forall X \in \{X\in\mathbb{S}^{n\times n},\mathcal{A}(X)=b\}, \quad \scal{M}{X}\leq \mu.
        \label{eq:separation_1}
      \end{gather}
    \end{subequations}
    Equation \eqref{eq:separation_2} is equivalent to $ M\succeq 0 \mbox{ and }\mu\leq 0.$
    %Because
    %\begin{equation*}
     % \{X\in \mathbb{S}^{n\times n},\mathcal{A}(X)=b\} = X_0 + \mathrm{Ker}(\mathcal{A}) = X_0 + \left(\mathrm{Range}(\mathcal{A}^*)\right)^\perp,
     % \end{equation*}
    And if we fix $V\in\mathcal{M}_p$, we can see that
    Equation \eqref{eq:separation_1} is equivalent to
    \begin{equation*}
      M \in\mathrm{Range}(\mathcal{A}^*)\mbox{ and } \scal{M}{VV^T} \leq \mu.
    \end{equation*}
    In particular,
    $ \scal{M}{VV^T} \leq \mu \leq 0.$
    As $M\succeq 0$, this means $MV=0$. Denoting $g\in\R^m$ a vector such that $M=\mathcal{A}^*(g)$, we have $ \mathcal{A}^*(g) V = 0$. Therefore, for any $\dot V\in\mathbb{R}^{n\times p}$,
    \begin{equation*}
      \scal{\mathcal{A}(V\dot V^T+\dot VV^T)}{g}
      =\scal{V\dot V^T+\dot VV^T}{M}
      = 2 \scal{\dot V}{MV}
      = 0,
    \end{equation*}
    which contradicts the assumption that $(\mathcal{A},b)$ is $p$-regular.

    \subsection{Proof of Proposition \ref{prop:strict_complementary_slackness}\label{app:strict_complementary_slackness}}

    We assume that strict complementary slackness holds, but $X_0$ is not the unique solution of Problem \eqref{eq:SDP}, and we show that $X_0$ is not an extremal point of $\mathcal{C}$.

    Let $X_0'$ be another solution. As $(C_1,g_1)$ is dual optimal, $X_0'$ and $C_1$ satisfy the complementary slackness condition:
    \begin{equation*}
      C_1 X_0'=0,
    \end{equation*}
    that is, $\mathrm{Range}(X_0')\subset \mathrm{Ker}(C_1)=\mathrm{Range}(X_0)$ (the last equality is because $C_1X_0=0$ and $\mathrm{rank}(X_0)+\mathrm{rank}(C_1)=n$).

    This inclusion and the fact that $X_0\succeq 0$ together imply that $X_0 + \epsilon(X_0'-X_0)\succeq 0$ for all $\epsilon\in\mathbb{R}$ close enough to $0$. And since $X_0,X_0'$ are both feasible points of Problem \eqref{eq:SDP},
    \begin{equation*}
      \mathcal{A}(X_0+\epsilon(X_0'-X_0))=b
    \end{equation*}
    for any $\epsilon\in\mathbb{R}$. Thus, $X_0+\epsilon(X_0'-X_0)$ is in the feasible set $\mathcal{C}$ of Problem \eqref{eq:SDP} for any $\epsilon$ close enough to $0$, and $X_0$ is not extremal.

\subsection{Proof of Proposition \ref{prop:first_order}\label{app:first_order}}

  From \citep*[Eq. 7]{boumal_voroninski_bandeira}, the gradient of the cost function of Problem \eqref{eq:P} at $V$ is $
    2 \mathrm{Proj}_V(CV),$
  where $\mathrm{Proj}_V:\R^{n\times p}\to T_V\mathcal{M}_p$ is the orthogonal projection onto $T_V\mathcal{M}_p$.
  Consequently, $V$ is a first-order critical point if and only if
  \begin{align*}
    CV \in (T_V\mathcal{M}_p)^\perp
       & = \{ \dot{V}\in\R^{n\times p}, \mathcal{A}(\dot{V}V^T+V\dot{V}^T)=0 \}^\perp\\
       & = \{ \dot{V}\in\R^{n\times p}, \forall g_2\in\R^m, \scal{\dot{V}V^T+V\dot{V}^T}{\mathcal{A}^*(g_2)}=0 \}^\perp\\
       & = \left(\{ \mathcal{A}^*(g_2)V, g_2\in\R^m\}^\perp \right)^\perp \\
    &= \{ \mathcal{A}^*(g_2)V, g_2\in\R^m\}.
  \end{align*}
  Now, $CV=\mathcal{A}^*(g_2)V$ for some $g_2\in\R^m$ if and only if $C = C_2 + \mathcal{A}^*(g_2)$, for some $g_2\in\R^m,C_2\in\mathbb{S}^{n\times n}$ such that $C_2 V = 0$.

  To show that, when it exists, the pair $(g_2,C_2)$ is unique, we assume that there exists another pair $(g'_2,C'_2)$ satisfying the same conditions. Then $\mathcal{A}^*(g_2-g_2')V = (C_2'-C_2)V=0$. The same argument as at the end of Appendix \ref{app:KKT_SDP} shows that $g_2-g_2'=0$. Therefore, $g_2=g_2'$ and $C_2=C_2'$.

  \subsection{Proof of Proposition \ref{prop:second_order}\label{app:second_order}}

  For any $\dot{V}\in T_V\mathcal{M}_p$, from \citep*[Eq. 10]{boumal_voroninski_bandeira},
  \begin{equation*}
    \mathrm{Hess}f_C(V) \cdot (\dot{V},\dot{V})
    = 2 \scal{S\dot{V}}{\dot{V}},
  \end{equation*}
  where $S=C-\mathcal{A}^*(\mu)$ for some $\mu\in\R^m$ such that $2SV=\mathrm{grad}f_C(V)=0$.

  From the uniqueness of $(C_2,g_2)$, we have $\mu=g_2$ and $S=C_2$.

\section{Discussion on Definition \ref{def:min_sec}\label{app:min_sec}}

Let for the time being $X_0\in\mathbb{S}^{n\times n},V\in \mathbb{R}^{n\times p}$ be fixed, such that $\mathrm{rank}(X_0)=r$. We assume Properties \ref{def:min_sec1} and \ref{def:min_sec2} of Definition \ref{def:min_sec} are true:
\begin{equation*}
  \mathrm{rank}(V)=p\quad\mbox{and}\quad
  \mathrm{Range}(X_0) \cap \mathrm{Range}(V)=\{0\}.
\end{equation*}
We discuss when Property \ref{def:min_sec3} holds. This property is equivalent to
\begin{align}
  T_V\mathcal{M}_p \cap \{\dot V,\mathrm{Range}(\dot V) \subset \mathrm{Range}(X_0)
  & + \mathrm{Range}(V)\} \nonumber\\
  & = \{VA,A\in\Anti(p)\}. \label{eq:min_sec3_equiv}
\end{align}
The vector spaces $\{\dot V,\mathrm{Range}(\dot V) \subset \mathrm{Range}(X_0) + \mathrm{Range}(V)\}$ and $T_V\mathcal{M}_p$ contain $\{VA,A\in\Anti(p)\}$, and respectively have dimensions
\begin{equation*}
  p\times \dim\left(\mathrm{Range}(X_0) + \mathrm{Range}(V)\right)
  =p(p+r)
\end{equation*}
and  $np-m$.

% Conversely, for any $(np-m)$-dimensional subspace $\mathcal{E}$ of $\mathbb{R}^{n\times p}$ containing $\{VA,A\in\Anti(p)\}$, one can show that there exists $(\mathcal{A},b)$ such that
%$T_V\mathcal{M}_p = \mathcal{E}$.

Consequently,
$\left(T_V\mathcal{M}_p\cap \{\dot V,\mathrm{Range}(\dot V) \subset \mathrm{Range}(X_0) + \mathrm{Range}(V)\}\right)^\perp$
is a subset of $\{VA,A\in\Anti(p)\}^\perp$, with dimension at most
\begin{align*}
  \min&\left( \dim\{VA,A\in\Anti(p)\}^\perp,
   \right.\\
      &\hskip 1cm \left. \dim(T_V\mathcal{M}_p)^\perp
        + \dim \{\dot V,\mathrm{Range}(\dot V) \subset \mathrm{Range}(X_0) + \mathrm{Range}(V)\}^\perp \right) \\
  & = \min\left(np-\frac{p(p-1)}{2}, m + np - p(p+r)\right) \\
  & = \dim\{VA,A\in\Anti(p)\}^\perp + \min\left(0, m  - \frac{p(p+1)}{2} - pr \right).
\end{align*}
Therefore, if $\frac{p(p+1)}{2}+pr>m$, $\left(T_V\mathcal{M}_p\cap \{\dot V,\mathrm{Range}(\dot V) \subset \mathrm{Range}(X_0) + \mathrm{Range}(V)\}\right)^\perp$ is a strict subset of $\{VA,A\in\Anti(p)\}^\perp$, Equation \eqref{eq:min_sec3_equiv} does not hold and $\mathcal{M}_p$ cannot be $X_0$-minimally secant at $V$. On the contrary, if $\frac{p(p+1)}{2}+pr\leq m$, the above upper bound on the dimension is exactly equal to the dimension if $T_V\mathcal{M}_p$ is ``generic enough'', hence $\mathcal{M}_p$ is $X_0$-minimally secant at $V$.

Consequently, we expect the main assumption in Theorem \ref{theorem:main} (the existence of $X_0,V$ such that $\mathcal{M}_p$ is $X_0$-minimally secant at $V$, in a setting where $\frac{p(p+1)}{2}+pr\leq m$) to hold for almost all $(\mathcal{A},b)$.

We however note that, although rare, there are pairs $(\mathcal{A},b)$ for which $X_0,V$ do not exist, hence the hypothesis cannot be trivially removed. An example is as follows.
\begin{example}
  We set $r=1, p=2$. Let $m\leq n$ be arbitrary. We consider
  \begin{equation*}
    \mathcal{A}:X\in\mathbb{S}^{n\times n}
    \to (X_{1,1},X_{1,n-m+2},\dots,X_{1,n})\in\mathbb{R}^m,
  \end{equation*}
  and $b=(1,0,\dots,0)$. One can check that $(\mathcal{A},b)$ is $2$-regular.

  The rank-$1$ elements of $\mathcal{C}$ are exactly the matrices $X_0$ of the form
  \begin{equation*}
    X_0 = \left(\begin{smallmatrix}
        \begin{smallmatrix}1&u^T\\u&uu^T\end{smallmatrix} & 0_{n-m+1,m-1} \\
        0_{m-1,n-m+1}&0_{m-1,m-1}
      \end{smallmatrix}\right)
    \quad\mbox{with }u\in\R^{(m-1)\times 1},
  \end{equation*}
  and $\mathcal{M}_2$ contains all matrices of the form $V=WX$, with $X\in O(2)$ and
  \begin{align*}
    W = 
    \left(\begin{smallmatrix} 1\,\, 0 \\ w_1\,w_2 \\ 0_{m-1,2}\end{smallmatrix}\right),
    \quad\mbox{with }w_1,w_2\in\R^{(n-m)\times 1}.
  \end{align*}
  For any rank-$1$ $X_0$ in $\mathcal{C}$ and $V\in\mathcal{M}_2$, using the above notations, one can check that
  \begin{equation*}
    \dot V =
    \left(\begin{smallmatrix}0&0\\\begin{smallmatrix}w_1-u\\\\0\end{smallmatrix}&\vdots\\\vdots&\vdots\\0&0\end{smallmatrix}\right)X
  \end{equation*}
  is in $T_V\mathcal{M}_2$, and $\mathrm{Range}(\dot V)\subset\mathrm{Range}(X_0)+\mathrm{Range}(V)$. Nevertheless, $\dot V\ne VA$, for all $A\in\Anti(2)$, so Property \ref{def:min_sec3} of Definition \ref{def:min_sec} does not hold (unless $w_1=u$, in which case Property \ref{def:min_sec2} does not hold), even if $m\geq \frac{p(p+1)}{2}+pr = 5$.
\end{example}

\section{Bad critical points can exist for $m<\frac{p(p+1)}{2}+pr$\label{app:non_tight}}

As announced in Remark \ref{rem:non_tight}, we provide an example where the conclusions of Theorem \ref{theorem:main} are true, but the assumption
\begin{equation*}
\frac{p(p+1)}{2}+pr \leq m
\end{equation*}
is not.

We set $r=p=2$ and $m=n=6$, $\mathcal{A}=\mathrm{diag}$ and $b=1_{6,1}$. To show that the conclusions of Theorem \ref{theorem:main} are true, it suffices to exhibit a matrix $C$ satisfying the conditions in Lemma \ref{lem:one_bad_C}.

  We set (this precise choice was suggested by numerical experiments)
  \begin{gather*}
    V = \left(\begin{smallmatrix}
      0&1\\
      0&1\\
      \frac{2}{\sqrt{5}}&\frac{1}{\sqrt{5}}\\
      \frac{1}{\sqrt{5}}&\frac{2}{\sqrt{5}}\\
      \frac{1}{\sqrt{5}}&\frac{2}{\sqrt{5}}\\
      0&1
    \end{smallmatrix}\right),\quad
    U_0 = \left(\begin{smallmatrix}
      \frac{2}{\sqrt{5}}&\frac{1}{\sqrt{5}}\\
      -1&0\\
      1&0\\
      -1&0\\
      \frac{2}{\sqrt{5}}&\frac{1}{\sqrt{5}}\\
      0&1
    \end{smallmatrix}\right),\quad
    X_0 =U_0U_0^T.
  \end{gather*}
  We also define
  \begin{gather*}
    g_1 = \begin{pmatrix}
      -\sqrt{5} &
      -2 + \frac{3}{\sqrt{5}} &
      -1 &
      -2 &
      0 &
      1
    \end{pmatrix}^T,\\
  g_2 = \begin{pmatrix}0&0&0&0&0&0\end{pmatrix}^T,\\
  C=C_2 = (G^{-1})^T
  \begin{pmatrix}
    0_{2,2}&0_{2,4}\\
      0_{4,2}&\begin{pmatrix}U_0&e_1&e_2\end{pmatrix}^T\mathrm{Diag}(g_1)\begin{pmatrix}U_0&e_1&e_2\end{pmatrix}
    \end{pmatrix}G^{-1} \\
    \hskip 4cm + 20 (G^{-1})^T
    \begin{pmatrix}
      0_{4,4}&0_{4,2}\\
      0_{2,4}&I_2
    \end{pmatrix}
    G^{-1},\\
    C_1=C-\mathrm{Diag}(g_1),
\end{gather*}
  where $e_1,e_2$ are the first two vectors of the canonical basis of $\R^{6\times 1}$, and $G=\begin{pmatrix}V&U_0&e_1&e_2\end{pmatrix}\in\R^{6\times 6}$ is the horizontal concatenation of $V,U_0,e_1,e_2$.

  With this choice, Properties \eqref{eq:12_cond1} and \eqref{eq:34_cond1} are true. We observe that
  \begin{equation*}
    C_2 V = C_2 G\begin{pmatrix}I_2\\0_{4,2}\end{pmatrix}=0_{6,2},
  \end{equation*}
  hence Property \eqref{eq:34_cond2} is also valid.
  A computation shows that
  \begin{align*}
    C_1 & = (G^{-1})^T
    \begin{pmatrix}
      \frac{6}{5}&\frac{6}{5} & 0&0 & 0&0\\
      \frac{6}{5}&\frac{14}{5}+\frac{2}{\sqrt{5}}&0&0&\sqrt{5}&2-\frac{3}{\sqrt{5}} \\
      0&0&0&0&0&0\\
      0&0&0&0&0&0\\
      0&\sqrt{5}&0&0&20&0\\
      0&2-\frac{3}{\sqrt{5}}&0&0&0&20
    \end{pmatrix} G^{-1}.
  \end{align*}
  From this expression, we see that Properties \eqref{eq:12_cond2} and \eqref{eq:12_cond4} are valid, as well as Property \eqref{eq:12_cond3} because
  \begin{equation*}
    C_1 U_0 = C_1 G\begin{pmatrix}0_{2,2}\\I_2\\0_{2,2}\end{pmatrix} = 0_{6,2}
    \qquad\Rightarrow\qquad C_1 X_0 = 0_{6,6}.
  \end{equation*}

  Finally, we consider Property \eqref{eq:34_cond3}. Let us define the bilinear form
  \begin{equation*}
  \begin{array}{rccc}
    q:& T_V\mathcal{M}_2 \times T_V\mathcal{M}_2&\to&\R\\
     & (\dot{V}_1,\dot{V}_2)&\to&\scal{C_2}{\dot{V}_1\dot{V}_2^T}.
  \end{array}
  \end{equation*}
  It contains $V\left(\begin{smallmatrix}0&1\\1&0\end{smallmatrix}\right)$ in its kernel (since $C_2V=0$). We can numerically compute the matrix associated to $q$ in an orthonormal basis of the $6$-dimensional space $T_V\mathcal{M}_2$ and check that it has $5$ strictly positive eigenvalues. Therefore, Property \eqref{eq:34_cond3} is also true.

  To summarize, Properties \eqref{eq:12_cond1} to \eqref{eq:34_cond3} are all true. The matrix $C$ therefore satisfies the properties required in Lemma \ref{lem:one_bad_C}.

  \section{Auxiliary results for the proof of Theorem \ref{theorem:bvb_improved}}

\subsection{Proof of Proposition \ref{prop:mathcalE}\label{app:mathcalE}}

Let $(V,C_2)$ belong to $\mathcal{E}$. We are going to exhibit a parametrization of $\mathcal{E}$ around $(V,C_2)$. Because $V$ has rank $p$, there exist a neighborhood $\mathcal{V}$ of $V$ in $\R^{n\times p}$ and a smooth map $W\in\mathcal{V}\to U_W\in\R^{n\times (n-p)}$ such that, for any $W\in\mathcal{V}$, the columns of $U_W$ form an orthonormal basis of $\mathrm{Range}(W)^\perp$. We define
\begin{equation*}
  \begin{array}{cccc}
    \psi : & (\mathcal{M}_p^{full}\cap \mathcal{V}) \times \mathbb{S}^{n\times n}
    & \to & \mathbb{S}^{p\times p} \times \mathbb{R}^{(n-p)\times p} \\
           & (W,D) & \to & (W^TDW,U_W^T DW).
  \end{array}
\end{equation*}
% \begin{equation*}
%   \begin{array}{cccc}
%     \psi : & (\mathcal{M}_p^{full}\cap \mathcal{V}) \times \mathbb{S}^{(n-p)\times (n-p)}
%     & \to & \mathcal{M}_p^{full} \times \mathbb{S}^{n\times n} \\
%            & (W,R) & \to & (W,U_W R U_W^T).
%   \end{array}
% \end{equation*}
This function is smooth. At any point $(W,D)$, its differential with respect to $D$ is surjective: for any $(A,B)\in\mathbb{S}^{p\times p}\times \mathbb{R}^{(n-p)\times p}$, one can check that $d_D\psi(W,D)\cdot \dot D = (A,B)$ if one sets
\begin{equation*}
  \dot D = (\begin{smallmatrix}W&U_W\end{smallmatrix})^{-1 T}\left(
    \begin{smallmatrix}
      A & B^T \\ B & 0_{n-p,n-p}
    \end{smallmatrix}\right) (\begin{smallmatrix}W&U_W\end{smallmatrix})^{-1}.
\end{equation*}
Therefore, from \citep*[Proposition 3.3.3]{absil}, $\psi^{-1}(0_{p,p},0_{n-p,p})$ is a submanifold of $\mathcal{M}_p^{full}\times\mathbb{S}^{n\times n}$, with dimension
\begin{equation*}
  \dim(\mathcal{M}_p^{full}\times \mathbb{S}^{n\times n})
  - \dim(\mathbb{S}^{p\times p}\times\mathbb{R}^{(n-p)\times p})
  = np-m + \frac{(n-p)(n-p+1)}{2}.
\end{equation*}

For any $(W,D)\in (\mathcal{M}_p^{full}\cap \mathcal{V}) \times \mathbb{S}^{n\times n}$, the following equivalences are true:
\begin{align*}
  \Big(\psi(W,D) = (0_{p,p},0_{n-p,p})\Big)
  &\quad \iff\quad \Big((\begin{smallmatrix}W&U_W\end{smallmatrix})^T DW = 0_{n,p}\Big) \\
  &\quad \iff\quad \Big( DW = 0_{n,p}\Big).
\end{align*}
Consequently, $\psi^{-1}(0_{p,p},0_{n-p,p})$ and $\mathcal{E}$ coincide in a neighborhood of $(V,C_2)$, which implies that $\mathcal{E}$ is also an $(np-m+\frac{(n-p)(n-p+1)}{2})$-dimensional manifold.

Its tangent space at $(V,C_2)$ is
\begin{align*}
  \mathrm{Ker}(d\psi(V,C_2))
  & = \{(\dot V,\dot C_2)\in T_V\mathcal{M}_p \times \mathbb{S}^{n\times n},
    \dot V^T C_2 V + V^T \dot C_2 V + V^T C_2 \dot V = 0_{p,p} \\
  & \hskip 2cm \mbox{and }
    (d U_V\cdot \dot V) C_2 V + U_V^T \dot C_2 V + U_V^T C_2 \dot V = 0_{n-p,p} \} \\
  & \overset{(C_2V=0)}= \{(\dot V,\dot C_2)\in T_V\mathcal{M}_p \times \mathbb{S}^{n\times n},
    V^T \dot C_2 V + V^T C_2 \dot V = 0_{p,p} \\
  & \hskip 2cm \mbox{and }
    U_V^T \dot C_2 V + U_V^T C_2 \dot V = 0_{n-p,p} \} \\
  & = \{(\dot V,\dot C_2)\in T_V\mathcal{M}_p \times \mathbb{S}^{n\times n},
    (\begin{smallmatrix}V & U_V\end{smallmatrix})^T (\dot C_2 V + C_2\dot V) = 0_{n,p}\} \\
  & = \{(\dot V,\dot C_2)\in T_V\mathcal{M}_p \times \mathbb{S}^{n\times n},
    \dot C_2 V + C_2\dot V = 0_{n,p}\}.
\end{align*}

\subsection{Proof of Proposition \ref{prop:face_Splus}\label{app:face_Splus}}

Let $U_V\in\R^{n\times p},U_V^{\perp}\in\R^{n\times (n-p)}$ be matrices whose columns respectively form an orthonormal basis of $\mathrm{Range}(V)$ and of $\mathrm{Ker}(VV^T)=\mathrm{Range}(V)^{\perp}$. Let $G\in\mathbb{R}^{p\times p}$ be the unique matrix such that $U_V = V G$.

Because $\mathcal{C}\subset \{X\in\mathbb{S}^{n\times n},X\succeq 0\}$, the face of $\mathcal{C}$ containing $VV^T$ is a subset of the face of $\{X\in\mathbb{S}^{n\times n},X \succeq 0\}$ containing $VV^T$, which is, from \citep*[Section 2.4]{laurent_rendl},
\begin{align*}
  \{X\in\mathbb{S}^{n\times n}, X\succeq 0,\mathrm{Ker}(VV^T)\subset \mathrm{Ker}(X)\}.
\end{align*}
Therefore, $\mathrm{Ker}(VV^T) \subset \mathrm{Ker}(X_{Face})$, which implies
\begin{align*}
  & X_{Face} U_V^{\perp} = 0_{n,n-p}, \\
  \Rightarrow\quad
  & \left(\begin{smallmatrix} U_V
      & U_V^{\perp} \end{smallmatrix}\right)^T
        X_{Face} \left(\begin{smallmatrix} U_V
            & U_V^{\perp} \end{smallmatrix}\right)
              = \left(\begin{smallmatrix} R & 0_{p,n-p} \\ 0_{n-p,p} & 0_{n-p,n-p}
                \end{smallmatrix}\right)\mbox{ for some }R\in\mathbb{S}^{p\times p}, \\
  \Rightarrow\quad
  & X_{Face}
    = \left(\begin{smallmatrix} U_V
        & U_V^{\perp} \end{smallmatrix}\right)
          \left(\begin{smallmatrix} R & 0_{p,n-p} \\ 0_{n-p,p} & 0_{n-p,n-p}
            \end{smallmatrix}\right) \left(\begin{smallmatrix} U_V
              & U_V^{\perp} \end{smallmatrix}\right)^T\mbox{ for some }R\in\mathbb{S}^{p\times p} \\
  % \Rightarrow\quad
  % & X_{Face} = U_V R U_V^T \mbox{ for some }R\in\mathbb{S}^{p\times p}, \\
  \Rightarrow\quad
  & X_{Face} = V GRG^T V^T \mbox{ for some }R\in\mathbb{S}^{p\times p}, \\
  \Rightarrow\quad
  & X_{Face} = V T V^T \mbox{ for some }T\in\mathbb{S}^{p\times p}.
\end{align*}

\section{Auxiliary results for the proof of Theorem \ref{theorem:main}}
  
  \subsection{Proof of Proposition \ref{prop:argmin_close_to_X0}\label{app:argmin_close_to_X0}}
  
  It suffices to show the following property:
    \begin{equation}\label{eq:cost_smaller_on_B}
      \mbox{for all }k\mbox{ large enough,}
      \quad\forall X\in \mathcal{C}-B(X_0,\epsilon),
      \quad \scal{C'_k}{X_0} < \scal{C'_k}{X}.
    \end{equation}
    Indeed, in this case, for $k$ large enough, any minimizer of $\scal{C'_k}{.}$ on the compact set $\mathcal{C}\cap \overline{B}(X_0,\epsilon)$ (there is at least one) is a minimizer of $\scal{C'_k}{.}$ on $\mathcal{C}$, and every minimizer of $\scal{C'_k}{.}$ on $\mathcal{C}$ is in $B(X_0,\epsilon)$.

    We assume, by contradiction, that Property \eqref{eq:cost_smaller_on_B} is not true. Up to replacing $(C'_k)_{k\in\mathbb{N}}$ by a subsequence, we can assume that, for any $k\in\mathbb{N}$,
    \begin{equation}\label{eq:def_X_prime_k}
      \exists X'_k \in\mathcal{C}-B(X_0,\epsilon), \quad \scal{C'_k}{X_0} \geq \scal{C'_k}{X'_k}.
    \end{equation}
    For any $k$, let $X'_k$ be such a matrix.

    By compactness, we can assume that $\left((X'_k-X_0)/||X'_k-X_0||\right)_{k\in\mathbb{N}}$
    converges to some unit-normed limit $Z\in\mathbb{S}^{n\times n}$. From Equation \eqref{eq:def_X_prime_k} and because $(C'_k)_{k\in\mathbb{N}}$ converges to $C$, $ \scal{C}{Z} \leq 0. $
    Equivalently,
    \begin{equation}\label{eq:X0_not_unique_min}
      \scal{C}{X_0+\epsilon Z} \leq \scal{C}{X_0}.
    \end{equation}
    Observe that $X_0+\epsilon Z$ belongs to $\mathcal{C}$: it is the limit of the sequence
    \begin{equation*}
      \left(\left(1-\frac{\epsilon}{||X'_k-X_0||}\right)X_0 + \frac{\epsilon}{||X'_k-X_0||}X'_k \right)_{k\in\mathbb{N}^*}.
    \end{equation*}
    Each element of this sequence belongs to $\mathcal{C}$ ($X_0$ and $X'_k$ do, and $\mathcal{C}$ is convex), and $\mathcal{C}$ is closed, so the limit also belongs to $\mathcal{C}$. Consequently, Equation \eqref{eq:X0_not_unique_min} contradicts the fact that $X_0$ is the unique minimizer of $\scal{C}{.}$ on $\mathcal{C}$.

  \subsection{Proof of Lemma \ref{lem:convergence_DC_hg}\label{app:convergence_DC_hg}}
  
  % \begin{lem*}[Lemma \ref{lem:convergence_DC_hg}]
  %   When $k$ goes to infinity,
  %   \begin{equation*}
  %     D_k \to C_1\quad\mbox{and}\quad h_k\to g_1.
  %   \end{equation*}
  % \end{lem*}
%  \begin{proof}[Proof of Lemma \ref{lem:convergence_DC_hg}]
    As $D_k + \mathcal{A}^*(h_k)=C'_k \overset{k\to+\infty}{\to} C=C_1+\mathcal{A}^*(g_1)$,
    \begin{equation}\label{eq:convergence_DChg}
      D_k-C_1 + \mathcal{A}^*(h_k-g_1) \overset{k\to+\infty}{\to} 0.
    \end{equation}
    In particular, if $h_k\overset{k\to+\infty}{\to}g_1$, then $D_k\overset{k\to+\infty}{\to}C_1$, so we only have to show that $(h_k)_{k\in\mathbb{N}}$ converges to $g_1$.

    By contradiction, we assume that $h_k \not \to g_1$.
    % Combining this hypothesis with a compactness argument, we will construct a non-zero $g \in\R^m$ such that
    % \begin{equation}\label{eq:h_infty_prop}
    %   \mathcal{A}^*(g)X_0 = 0,
    % \end{equation}
    % which is impossible from Proposition \ref{prop:Astar_g_X0}.    
    Up to replacing $(h_k)_{k\in\mathbb{N}}$ by a subsequence, we can assume that $(||h_k-g_1||)_{k\in\mathbb{N}}$ is lower bounded by a positive constant and that $((h_k-g_1)/||h_k-g_1||)_{k\in\mathbb{N}}$ converges to some non-zero limit $g$.
    %. Let $g\in\R^m$ be the limit, and establish Property \eqref{eq:h_infty_prop}.

    From Equation \eqref{eq:convergence_DChg},
    \begin{equation}\label{eq:convergence_CD_h_infty}
      \frac{C_1-D_k}{||h_k-g_1||} \overset{k\to+\infty}{\to} \mathcal{A}^*(g).
    \end{equation}
    From Proposition \ref{prop:argmin_close_to_X0}, $(X'_k)_{k\in\mathbb{N}}$ converges to $X_0$, so $
      C_1 X'_k \overset{k\to+\infty}{\to} C_1 X_0 = 0,$ 
    and because $(||h_k-g_1||)_{k\in\mathbb{N}}$ is bounded away from zero, this implies
    \begin{equation*}
      \frac{C_1X'_k}{||h_k-g_1||} \overset{k\to+\infty}{\to} 0.
    \end{equation*}
    Recalling that, from the definition of $D_k$, $D_kX'_k=0$ for all $k$, Equation \eqref{eq:convergence_CD_h_infty} yields:
    \begin{equation*}
      \mathcal{A}^*(g)X_0 = \underset{k\to+\infty}{\lim}\left(\frac{C_1-D_k}{||h_k-g_1||}\right)X'_k = 0.
    \end{equation*}

    Therefore, we also have $\mathcal{A}^*(g)V_0$, if we fix $V_0\in\mathbb{R}^{n\times p}$ such that $X_0=V_0V_0^T$ (it is possible, as $X_0\succeq 0$ and $\mathrm{rank}(X_0)=r\leq p$). The matrix $V_0$ is in $\mathcal{M}_p$. Applying the same argument as at the end of Appendix \ref{app:KKT_SDP}, we reach a contradiction.
    
%  \end{proof}

    \subsection{Proof of Proposition \ref{prop:second_order_stability_general}\label{app:second_order_stability_general}}

    If we compose $f$ with a diffeomorphism along the second coordinate, we can assume that $\mathcal{M}$ is an open subset of $\R^d$, for some integer $d$ (from Proposition \ref{prop:transfer_critical_points}, the critical points of the composition of a function and a diffeomorphism are exactly the image by the reciprocal diffeomorphism of the critical points of the function).

    Let $\mathcal{V}$ be a neighborhood of $v$ in $\mathcal{M}$.
    We define
    \begin{equation*}
      \begin{array}{cccc}
        \chi:& E \times \mathcal{V} &\to&\R^d \\
             & (c',v') & \to &\nabla(f(c',.))(v').
      \end{array}
    \end{equation*}
    This is a smooth map; it satisfies $\chi(c,v)=0$ (since $v$ is a critical point of $f(c,.)$) and its differential at $(c,v)$ along the second coordinate is invertible (it is $\mathrm{Hess}(f(c,.))(v)$, which is positive definite by assumption). From the implicit function theorem, there exist a neighborhood $\mathcal{E}$ of $c$ in $E$, and a smooth function $\delta:\mathcal{E}\to \mathcal{V}$ such that $\delta(c) = v$ and $\chi(c',\delta(c'))=0$ for any $c'\in\mathcal{E}$. We fix such $\mathcal{E},\delta$. Then, for any $c'\in\mathcal{E}$,
    \begin{equation*}
      \nabla(f(c',.))(\delta(c')) = \chi(c',\delta(c')) = 0.
    \end{equation*}
    Equivalently, $\delta(c')$ is a first-order critical point of $f(c',.)$. Additionally, the map $c'\to \mathrm{Hess}(f(c',.))(\delta(c'))$ is continuous ($f$ and $\delta$ are smooth), and
    \begin{equation*}
      \mathrm{Hess}(f(c,.))(\delta(c)) = \mathrm{Hess}(f(c,.))(v) \succ 0.
    \end{equation*}
    As a consequence, for any $c'\in\mathcal{E}$ close enough to $c$,
    \begin{equation*}
      \mathrm{Hess}(f(c',.))(\delta(c' )) \succ 0.
    \end{equation*}
    Therefore, for any $c'$ close enough to $c$, $\delta(c')$, which is an element of $\mathcal{V}$, is a second-order critical point of $f(c',.)$.    

  \subsection{Proof of Proposition \ref{prop:transfer_critical_points}\label{app:transfer_critical_points}}

  Let $v$ belong to $\mathcal{N}_1$. We have
    \begin{equation*}
      \nabla(f\circ \phi)(v) =
      (d\phi(v))^* \nabla f (\phi(v)).
    \end{equation*}
    As $d\phi(v)^*$ is injective (it is the adjoint of a surjective map), $v$ is a first-order critical point of $f\circ \phi$ if and only if $\phi(v)$ is a first-order critical point of $f$.

    In this case, the Hessians of $f$ and $f\circ\phi$ at $\phi(v)$ and $v$ are linked by the following relation:
    \begin{align*}
      \forall x_1,x_2\in T_v\mathcal{N}_1,\quad
      \mathrm{Hess}(f\circ\phi)(v)\cdot(x_1,x_2)
      = \mathrm{Hess} f(\phi(v))\cdot (d\phi(v)\cdot x_1,d\phi(v)\cdot x_2).
    \end{align*}
    As $d\phi(v)$ is surjective, $\mathrm{Hess}(f\circ \phi)(v)$ and $\mathrm{Hess}f(\phi(v))$ have the same rank, and $\mathrm{Hess}(f\circ \phi)(v)$ is positive semidefinite if and only if $\mathrm{Hess}f (\phi(v))$ is, meaning that $v$ is a second-order critical point of $f\circ \phi$ if and only if $\phi(v)$ is a second-order critical point of $f$.

    \subsection{Proof of Proposition \ref{prop:up_to_replacing}\label{app:up_to_replacing}}

 Let us consider for a moment an arbitrary invertible matrix $G\in\mathbb{R}^{n\times n}$, and define
\begin{gather*}
  \tilde{\mathcal{A}} : X\in\mathbb{S}^{n\times n}
  \to \mathcal{A}(GXG^T) \in\mathbb{R}^m
  \quad\mbox{and}\quad
  \tilde b = b, \\
  \tilde X_0 = G^{-1} X_0 (G^T)^{-1}
  \quad\mbox{and}\quad
  \tilde V = G^{-1}V.
\end{gather*}
We denote $\tilde{\mathcal{M}}_p$ the set of feasible points for Problem \eqref{eq:P} where $\mathcal{A}$ and $b$ have been replaced with $\tilde{\mathcal{A}}$ and $\tilde{b}$:
\begin{align*}
  \tilde{\mathcal{M}}_p
  & = \{W\in\mathbb{R}^{n\times p},\tilde{\mathcal{A}}(WW^T)=\tilde{b}\} \\
  & = \{ G^{-1}W, W\in\mathcal{M}_p\}.
\end{align*}
The pair $(\tilde{\mathcal{A}},\tilde b)$ is $p$-regular (because $(\mathcal{A},b)$ is) and one can check that $\tilde{\mathcal{M}}_p$ is $\tilde{X}_0$-minimally secant at $\tilde{V}$ (because $\mathcal{M}_p$ is $X_0$-minimally secant at $V$).

Now imagine that we can construct $\tilde{C}_1,\tilde{g}_1,\tilde{C}_2,\tilde{g}_2$ satisfying Conditions \eqref{eq:12_cond2} to \eqref{eq:34_cond3} (with $\tilde{X}_0,\tilde{V},\tilde{\mathcal{A}},\tilde{\mathcal{M}}_p$ in place of their non-tilde versions). Then, if we define
\begin{equation*}
  C_1 = (G^T)^{-1} \tilde{C}_1 G^{-1},\quad
  C_2 = (G^T)^{-1} \tilde{C}_2 G^{-1},\quad
  g_1 = \tilde{g}_1,\quad
  g_2 = \tilde{g}_2,
\end{equation*}
we see that these objects satisfy Conditions \eqref{eq:12_cond2} to \eqref{eq:34_cond3} (with the non-tilde versions this time). Therefore, if we are able to construct $\tilde{C}_1,\tilde{g}_1,\tilde{C}_2,\tilde{g}_2$ satisfying Conditions \eqref{eq:12_cond2} to \eqref{eq:34_cond3}, it proves the existence of $C_1,g_1,C_2,g_2$ satisfying these same conditions.

To conclude, it suffices to show that, if we properly define $G$, then
\begin{equation}
  \tilde X_0 = \left(\begin{smallmatrix}
      I_r&0_{r,n-r}\\
      0_{n-r,r}&0_{n-r,n-r}
    \end{smallmatrix}\right)
  \quad\mbox{and}\quad
  \tilde V = \left(\begin{smallmatrix}0_{r,p} \\ I_p \\ 0_{n-p-r,p} \end{smallmatrix}\right).
  \label{eq:tildeX_tildeV}
\end{equation}
Let $U_0\in\mathbb{R}^{n\times r}$ be such that $X_0=U_0U_0^T$ (it exists: $X_0$ is semidefinite positive and has rank $r$). We define
\begin{equation*}
  G = \begin{pmatrix}U_0 & V & W \end{pmatrix} \in\mathbb{R}^{n\times n},
\end{equation*}
where $W\in\mathbb{R}^{n\times (n-r-p)}$ is any matrix that makes $G$ invertible (it exists, as the columns of $U_0$ and $V$ are linearly independent, from Properties \ref{def:min_sec1} and \ref{def:min_sec2} of Definition \ref{def:min_sec}). Equation \eqref{eq:tildeX_tildeV} holds.

\subsection{Proof of Lemma \ref{lem:construction_g1}\label{app:construction_g1}}

We define $L:\mathbb{S}^{n\times n}\to \mathbb{R}^{r\times p}\times\mathbb{S}^{p\times p}$ the linear map such that, for any $R_1,R_2,G_1,G_2,G_4,G_5$,
\begin{equation*}
L\left(\left(\begin{smallmatrix}G_1&R_1&G_2
        \\ R_1^T&R_2&G_4 \\ G_2^T&G_4^T&G_5
      \end{smallmatrix}\right)\right) = (R_1,R_2).
\end{equation*}
Proving the lemma amounts to showing that $L\circ\mathcal{A}^*$ is surjective. Equivalently, it suffices to show that the dual map $\mathcal{A}\circ L^*$ is injective. Let $(R_1,R_2)$ be in its kernel:
\begin{equation*}
  \mathcal{A}
  \left(\left(\begin{smallmatrix}0&R_1/2&0
        \\ R_1^T/2&R_2&0 \\ 0&0&0
      \end{smallmatrix}\right)\right) = 0.
\end{equation*}
We recall that we have assumed, following Proposition \ref{prop:up_to_replacing},
\begin{equation*}
  X_0 = \left(\begin{smallmatrix}
      I_r&0_{r,n-r}\\
      0_{n-r,r}&0_{n-r,n-r}
    \end{smallmatrix}\right)
  \quad\mbox{and}\quad
  V = \left(\begin{smallmatrix}0_{r,p} \\ I_p \\ 0_{n-p-r,p} \end{smallmatrix}\right).
\end{equation*}
Therefore, if we set
\begin{equation*}
  \dot V = \left(\begin{smallmatrix}
      R_1/2\\R_2/2\\0_{n-p-r,p}
    \end{smallmatrix} \right),
\end{equation*}
we have $\mathrm{Range}(\dot V)\subset\mathrm{Range}(X_0)+\mathrm{Range}(V)$ and $\mathcal{A}(V\dot V^T+\dot V V^T)=0$, hence $\dot V\in T_V\mathcal{M}_p$. From Property \ref{def:min_sec3} of Definition \ref{def:min_sec}, there exists $A\in\Anti(p)$ such that $\dot V=VA$. As a consequence, $R_1=0_{r,p}$ and $R_2$ is both symmetric and antisymmetric, hence $R_2=0_{p,p}$.
This proves that $\mathrm{Ker}(\mathcal{A}\circ L^*)=\{(0_{r,p},0_{p,p})\}$, which is what we needed.

\subsection{Proof of Proposition \ref{prop:properties_still_hold}\label{app:properties_still_hold}}

Properties \eqref{eq:12_cond2}, \eqref{eq:12_cond3}, \eqref{eq:34_cond1} and \eqref{eq:34_cond2} are a direct consequence of Equation \eqref{eq:up_to_replacing} and of the fact that $g_1,C_1,C_2$ satisfy these same properties.

For Property \eqref{eq:12_cond4}, we have $\mathrm{rank}(C_1^{(mod)})\geq \mathrm{rank}(C_1)=n-\mathrm{rank}(X_0)$, since adding a semidefinite positive matrix to another one cannot decrease the rank. Additionally, as $C_1^{(mod)}X_0=0$ (Property \eqref{eq:12_cond3}), we also have $\mathrm{rank}(C_1^{(mod)})\leq n -\mathrm{rank}(X_0)$ and, therefore, $\mathrm{rank}(C_1^{(mod)})= n -\mathrm{rank}(X_0)$.

\subsection{Proof of Proposition \ref{prop:34_cond3_orth1}\label{app:34_cond3_orth1}}

Equation \eqref{eq:34_cond3_bis} naturally implies \eqref{eq:34_cond3_bis_orth}. Let us assume that Equation \eqref{eq:34_cond3_bis_orth} is true and show the converse.

Let $\dot V$ be in $T_V\mathcal{M}_p$. We must show that $\scal{C_2^{(mod)}}{\dot V \dot V^T}\geq 0$, with equality if and only if $\dot V = VA$ for some $A\in\Anti(p)$. We write
\begin{equation*}
  \dot V = \dot W + VA
  \mbox{ for some }\dot W\in\mathcal{E}_{\perp},A\in\Anti(p).
\end{equation*}
Using at the last line the fact that $C_2^{(mod)}V=0$ (Property \eqref{eq:34_cond2}), we see that
\begin{align*}
  \scal{C_2^{(mod)}}{\dot V \dot V^T}
  & = \scal{C_2^{(mod)}}{\dot W \dot W^T}
    + \scal{C_2^{(mod)}}{V A \dot W^T} \\
  & \hskip 1cm
    + \scal{C_2^{(mod)}}{\dot W A^T V^T}
    + \scal{C_2^{(mod)}}{V A A^T V^T} \\
  & = \scal{C_2^{(mod)}}{\dot W \dot W^T}
    + 2\scal{C_2^{(mod)}V}{\dot W A^T} 
    + \scal{C_2^{(mod)}V}{V A A^T} \\
  & = \scal{C_2^{(mod)}}{\dot W\dot W^T}.
\end{align*}
Therefore, from Equation \eqref{eq:34_cond3_bis_orth}, $\scal{C_2^{(mod)}}{\dot V \dot V^T} \geq 0$, with equality if and only if $\dot W=0$, that is $\dot V \in \{VA,A\in\Anti(p)\}$.

\subsection{Proof of Proposition \ref{prop:34_cond3_orth2}\label{app:34_cond3_orth2}}

First, we observe that for any $\dot V\in\mathcal{E}_{\perp}-\{0\}$,
\begin{equation*}
\scal{\left(\begin{smallmatrix}
      0_{r+p,r+p}&0_{r+p,n-r-p}\\
      0_{n-r-p,r+p}&I_{n-r-p}
    \end{smallmatrix}\right)}{\dot V \dot V^T} \geq 0
\end{equation*}
because it is the scalar product of two semidefinite positive matrices. It is zero if and only if the last $n-r-p$ rows of $\dot V$ are zero, that is
\begin{equation*}
\mathrm{Range}(\dot V) \subset \mathrm{Range}(X_0) + \mathrm{Range}(V).
\end{equation*}
Because $\mathcal{M}_p$ is $X_0$-minimally secant at $V$, this is possible only if $\dot V = VA$ for some $A\in\Anti(p)$, which contradicts the fact that $\dot V$ is in $\mathcal{E}_{\perp}-\{0\}$. Therefore,
\begin{equation*}
\scal{\left(\begin{smallmatrix}
      0_{r+p,r+p}&0_{r+p,n-r-p}\\
      0_{n-r-p,r+p}&I_{n-r-p}
    \end{smallmatrix}\right)}{\dot V \dot V^T} > 0.
\end{equation*}

We set $\mathcal{B}_{\perp}=\{\dot V\in\mathcal{E}_{\perp}, ||\dot V||_F=1\}$. From the previous remark and because $\mathcal{B}_{\perp}$ is compact, there exists $\epsilon>0$ such that
\begin{equation*}
  \forall \dot V\in\mathcal{B}_{\perp},\quad
  \scal{\left(\begin{smallmatrix}
      0_{r+p,r+p}&0_{r+p,n-r-p}\\
      0_{n-r-p,r+p}&I_{n-r-p}
    \end{smallmatrix}\right)}{\dot V \dot V^T} \geq \epsilon.
\end{equation*}
We define $\gamma= \inf_{\dot V \in\mathcal{B}_{\perp}}\scal{C_2}{\dot V\dot V^T}$. For any $t$ such that $\gamma + t\epsilon >0$, it holds:
\begin{equation*}
  \forall \dot V\in\mathcal{B}_{\perp},\quad
  \scal{C_2^{(mod)}}{\dot V \dot V^T} \geq \gamma + t\epsilon >0.
\end{equation*}
In this case, by homogeneity, Equation \eqref{eq:34_cond3_bis_orth} also holds.

\section{Proofs for Subsection \ref{ss:examples}}

  \subsection{Proof of Corollary \ref{cor:maxcut}\label{app:maxcut}}

  % \begin{cor*}[Corollary \ref{cor:maxcut}]
  %   If $p\in\mathbb{N}$ is such that $
  %   \frac{p(p+1)}{2}+p > n $
  %   then, for almost any cost matrix $C$, all second-order critical points of the Burer-Monteiro factorization of Problem \eqref{eq:SDP-MaxCut} are globally optimal.  
    
  %   On the other hand, for any $p$ such that $
  %   \frac{p(p+1)}{2} + p \leq n $
  %   the set of cost matrices admits a subset with non-zero Lebesgue measure on which
  %   \begin{itemize}
  %   \item Problem \eqref{eq:SDP-MaxCut} has a unique global minimizer, which has rank $1$;
  %   \item Its Burer-Monteiro factorization with rank $p$ has at least one non-optimal second-order critical point.
  %   \end{itemize}
  % \end{cor*}

  % \begin{proof}[Proof of Corollary \ref{cor:maxcut}]

  The first part of the corollary is a direct consequence of Theorem \ref{theorem:bvb_improved}, so we focus on the second one.
 Let $p,n$ be such that $p(p+1)/2+p\leq n$. It suffices to check the three hypotheses of Theorem \ref{theorem:main}. The first two are classical, setting
  \begin{equation*}
    X_0=U_0 U_0^T\qquad\mbox{with}\qquad U_0=\left(\begin{smallmatrix}1\\\vdots\\1\end{smallmatrix}\right) \in\R^{n\times 1}.
  \end{equation*}

  Let us construct $V\in\mathcal{M}_p$ such that $\mathcal{M}_p$ is $X_0$-minimally secant at $V$. In Definition \ref{def:min_sec}, the only delicate part is Property \ref{def:min_sec3}. We do not have a better method to check it that direct computation. Hence, we must choose $V$ as simple as possible, so that the equations defining $T_V\mathcal{M}_p$ and $\{\dot V\in\mathbb{R}^{n\times p},\mathrm{Range}(\dot V)\subset \mathrm{Range}(U_0) + \mathrm{Range}(V)\}$ are relatively easy to manipulate. This matrix must satisfy two constraints: it has to be in $\mathcal{M}_p$ (that is, all its rows must have norm $1$) and it must have at least $\frac{p(p+1)}{2}+p$ different lines (otherwise, one can check that the aforementioned equations are degenerate).

  The simplest matrix $V$ that satisfies these constraints is arguably the following one: for any $i\leq p$ and $j\in\{i+1,\dots,p\}$, we respectively set the $i$-th, $(p+i)$-th and $(2p+\phi(i,j))$-th lines of $V$ as
  \begin{equation}\label{eq:maxcut_def_V}
    V_{i,:}=e_i,\quad V_{p+i,:}=-e_i,\quad
    V_{2p+\phi(i,j),:}=\frac{e_i+e_j}{\sqrt{2}},
  \end{equation}
where $(e_1,\dots,e_p)$ is the canonical basis of $\mathbb{R}^{1\times p}$ and $\phi:\{i,j\mbox{ s.t. }1\leq i<j\leq p\} \to \left\{1,\dots,\frac{p(p-1)}{2}\right\}$ is an arbitrary bijection. For the last $n-\left(\frac{p(p+1)}{2}+p\right)$ lines, we choose any unit-normed elements of $\R^{1\times p}$.

  This definition ensures that $V$ has rank $p$ (it contains $I_p$ as a submatrix). Moreover, $\left(\begin{smallmatrix}U_0&V\end{smallmatrix}\right)$ has rank $p+1$ (its $p+1$ first lines form an invertible matrix). Therefore, Properties \ref{def:min_sec1} and \ref{def:min_sec2} of Definition \ref{def:min_sec} hold.

  We check Property \ref{def:min_sec3}. Let $\dot{V}\in T_V\mathcal{M}_p$ be such that
  \begin{equation*}
    \mathrm{Range}(\dot V)\subset \mathrm{Range}(U_0)+\mathrm{Range}(V).
  \end{equation*}
  Then there exists $(R,A)\in\R^{1\times p}\times \R^{p\times p}$ such that $\dot V = U_0 R + V A$. We fix such $R,A$, and show that $R=0_{1,p}$ and $A$ is antisymmetric.

  For any $i=1,\dots,p$, because $\dot V$ is in $T_V\mathcal{M}_p$,
  \begin{align*}
    &\left(\mathrm{diag}(V\dot V^T+\dot V V^T)_i = 0 \mbox{ and }
    \mathrm{diag}(V\dot V^T+\dot V V^T)_{p+i} = 0\right) \\
    \overset{\textrm{\eqref{eq:maxcut_def_V}}}{\iff}\quad
    &\left(\dot V_{i,i}=0 \mbox{ and }\dot V_{p+i,i}=0 \right) \\
    \iff\quad
    &\left( (U_0R+VA)_{i,i}=0\mbox{ and } (U_0R+VA)_{p+i,i}=0 \right) \\
    \overset{\textrm{\eqref{eq:maxcut_def_V}}}{\iff}\quad
    &\left(R_{1,i}+A_{i,i}=0 \mbox{ and } R_{1,i}-A_{i,i}=0 \right) \\
    \iff\quad
    &\left(R_{1,i}=0\mbox{ and }A_{i,i}=0\right).
  \end{align*}
  Consequently, $R=0_{1,p}$ and $\mathrm{diag}(A)=0$. Similarly, for any $1\leq i<j\leq p$,
  \begin{align*}
    &\left(\mathrm{diag}(V\dot V^T+\dot V V^T)_{2p+\phi(i,j)}=0\right) \\
    \overset{\textrm{\eqref{eq:maxcut_def_V}}}{\iff} \quad
    &\left( \dot V_{2p+\phi(i,j),i} + \dot V_{2p+\phi(i,j),j}=0 \right) \\
    \iff\quad
    &\left((U_0R+VA)_{2p+\phi(i,j),i} + (U_0R+VA)_{2p+\phi(i,j),j}=0\right) \\
    \overset{\textrm{\eqref{eq:maxcut_def_V}}}{\iff} \quad
    &\left( R_{1,i} + \frac{A_{i,i}+A_{i,j}}{\sqrt{2}}
      +R_{1,j}+\frac{A_{j,i}+A_{j,j}}{\sqrt{2}} = 0 \right) \\
    \iff\quad
    &\left(A_{i,j}+A_{j,i}=0\right).
  \end{align*}
  The matrix $A$ is therefore antisymmetric.

  \subsection{Proof of Corollary \ref{cor:orthogonal_cut}\label{app:orthogonal_cut}}

  The first part of the corollary is a direct consequence of Theorem \ref{theorem:bvb_improved}.

  Let us fix $S$ and $p\geq d$ such that
  $\frac{p(p+1)}{2} + pd \leq \frac{Sd(d+1)}{2}$, and prove the second part by checking that the hypotheses of Theorem \ref{theorem:main} hold true.
  We have already said that $(\mathcal{A},b)$ is $p$-regular, which is the second hypothesis.
  
  We set
  \begin{equation*}
    U_0 = \left(\begin{smallmatrix}
        I_d \\ \vdots \\ I_d
      \end{smallmatrix}\right) \in \R^{Sd\times d}
    \quad\mbox{and}\quad
    X_0 = U_0 U_0^T \in\mathbb{S}^{Sd\times Sd}.
  \end{equation*}
  Then $X_0$ is an extreme point of $\mathcal{C}$ with rank $d$: the first hypothesis holds.

  We show the third hypothesis by exhibiting $V$ such that $\mathcal{M}_p$ is $X_0$-minimally secant at $V$. We did not find a general construction that would be applicable for any value of $d$. Hence, we present separate constructions for the cases $d=1$, $d=2$ and $d=3$.

  The case $d=1$ has already been studied in Subsection \ref{app:maxcut}, so we consider the case where $d=2$.
    We define the following blocks:
    \begin{gather*}
      G_1 = \left(\begin{smallmatrix}
          1&0&0 \\ 0&1&0
        \end{smallmatrix}\right), \qquad
      G_2 = \left(\begin{smallmatrix}
          0&1&0 \\ 0&0&1
        \end{smallmatrix}\right), \qquad      
      G_3 = \left(\begin{smallmatrix}
          0&0&1 \\ \frac{1}{\sqrt{2}}&\frac{1}{\sqrt{2}}&0
        \end{smallmatrix}\right), \\
      G_4 = \left(\begin{smallmatrix}
          0&\frac{1}{\sqrt{2}}&\frac{1}{\sqrt{2}} \\
          \frac{1}{\sqrt{3}}&\frac{1}{\sqrt{3}}&-\frac{1}{\sqrt{3}}
        \end{smallmatrix}\right).
    \end{gather*}
    
    We distinguish depending on the congruency of $p$ modulo $3$.
    If $p\equiv 0[3]$, for any $q=1,\dots,p/3$, we set
    \begin{gather*}
      W_q^{(1)} = \left( \begin{smallmatrix}
          0_{2\times 3(q-1)}&G_1&0_{2\times (p-3q)}
        \end{smallmatrix} \right),\qquad
      W_q^{(2)} = \left( \begin{smallmatrix}
          0_{2\times 3(q-1)}&G_2&0_{2\times (p-3q)}
        \end{smallmatrix} \right),\\
      W_q^{(3)} = \left( \begin{smallmatrix}
          0_{2\times 3(q-1)}&G_3&0_{2\times (p-3q)}
        \end{smallmatrix} \right),\qquad
      W_q^{(4)} = \left( \begin{smallmatrix}
          0_{2\times 3(q-1)}&G_4&0_{2\times (p-3q)}
        \end{smallmatrix} \right).
    \end{gather*}
    For any $q,q'\in\{1,\dots,p/3\}$ such that $q<q'$, we set
    \begin{gather*}
      X_{q,q'}^{(1)} = \left( \begin{smallmatrix}
          0_{2\times 3(q-1)}&G_1&0_{2\times 3(q'-q-1)}&G_1&0_{2\times (p-3q')}
        \end{smallmatrix} \right) / \sqrt{2}, \\
      X_{q,q'}^{(2)} = \left( \begin{smallmatrix}
          0_{2\times 3(q-1)}&G_2&0_{2\times 3(q'-q-1)}&G_3&0_{2\times (p-3q')}
        \end{smallmatrix} \right) / \sqrt{2}, \\
      X_{q,q'}^{(3)} = \left( \begin{smallmatrix}
          0_{2\times 3(q-1)}&G_4&0_{2\times 3(q'-q-1)}&G_2&0_{2\times (p-3q')}
        \end{smallmatrix} \right) / \sqrt{2}.
    \end{gather*}
    We define
    \begin{equation*}
      V = \left(\begin{smallmatrix}
          W_1^{(1)}\\W_1^{(2)}\\
          \vdots \\
          X_{1,2}^{(1)}\\X_{1,2}^{(2)}\\ \vdots \\
          X_{p/3-1,p/3}^{(3)}\\
          V_{\frac{p^2+5p}{6}+1}\\\vdots\\V_S
        \end{smallmatrix}\right),
    \end{equation*}
    where $V_{\frac{p^2+5p}{6}+1},\dots,V_S$ are arbitrary elements of $\mathbb{R}^{2\times p}$ such that, for all $k$, $V_kV_k^T=I_2$.

    Let $\dot V\in T_V\mathcal{M}_p$ be such that
    \begin{equation}\label{eq:orth_cut_inclusion_ranges}
      \mathrm{Range}(\dot V)\subset \mathrm{Range}(X_0) + \mathrm{Range}(V)
      = \mathrm{Range}(U_0) + \mathrm{Range}(V).
    \end{equation}
    We show that $\dot V=VA$ for some $A\in\mathrm{Anti}(p)$. Equation \eqref{eq:orth_cut_inclusion_ranges} means that there exists $R\in\mathbb{R}^{2\times p},T\in\mathbb{R}^{p\times p}$ such that
    \begin{equation*}
      \dot V = U_0 R + V T.
    \end{equation*}
    We call $R_1,\dots,R_{p/3}$ the elements of $\mathbb{R}^{2\times 3}$ and $T_{1,1},\dots,T_{p/3,p/3}$ the ones of $\mathbb{R}^{3\times 3}$ such that
    \begin{equation*}
      R = \begin{pmatrix} R_1&\dots&R_{p/3}\end{pmatrix}
      \quad\mbox{and}\quad
      T = \left(\begin{smallmatrix} T_{1,1}&\dots&T_{1,p/3}\\
          \vdots&&\vdots \\
          T_{p/3,1}&\dots&T_{p/3,p/3}
        \end{smallmatrix}\right),
    \end{equation*}
    
    Because $\dot V$ belongs to $T_V\mathcal{M}_p$, the $k$-th $2\times 2$ diagonal block of $\left(\dot V V^T + V \dot V^T\right)$ is zero for any $k=1,\dots,S$, that is, the $k$-th block of $\dot V V^T$ is antisymmetric. Using the definitions of $U_0$ and $V$, this property, for $k\leq \frac{4}{3}p$, can be rewritten as
    \begin{equation*}
      R_qG_s^T + G_sT_{q,q}G_s^T
      \in\mathrm{Anti}(2),
      \quad\forall q\leq \frac{p}{3},s\leq 4.
      \label{eq:orth_cut_first_system}
    \end{equation*}
    For any fixed $q$, this property implies that $R_q,T_{q,q}$, which together contain $15$ coefficients, satisfy a linear system with $12$ equality constraints (there are four possible values for $s$, and each value yields $3$ constraints). It is tedious but not difficult to check that $(R_q,T_{q,q})$ is a solution of this linear system if and only if
    \begin{equation}\label{eq:orth_cut_R_Tqq}
      R_q=0_{2,3}\quad\mbox{and}\quad T_{q,q}\in\mathrm{Anti}(3).
    \end{equation}
    Now, considering the $k$-th diagonal block of $\dot V V^T$ for $k=\frac{4}{3}p+1,\dots, \frac{p^2+5p}{6}$, we see that, for any $1\leq q<q'\leq \frac{p}{3}$,
    \begin{align*}
      \frac{R_qG_1^T + R_{q'}G_1^T}{\sqrt{2}}
      + \frac{G_1 T_{q,q}G_1^T + G_1 T_{q,q'}G_1^T + G_1 T_{q',q} G_1^T + G_1 T_{q',q'} G_1^T}{2},\\
      \frac{R_qG_2^T + R_{q'}G_3^T}{\sqrt{2}}
      + \frac{G_2 T_{q,q}G_2^T + G_2 T_{q,q'}G_3^T + G_3 T_{q',q} G_2^T + G_3 T_{q',q'} G_3^T}{2},\\
      \frac{R_qG_4^T + R_{q'}G_2^T}{\sqrt{2}}
      + \frac{G_4 T_{q,q}G_4^T + G_4 T_{q,q'}G_2^T + G_2 T_{q',q} G_4^T + G_4 T_{q',q'} G_4^T}{2}
    \end{align*}
    are also antisymmetric. Taking Equation \eqref{eq:orth_cut_R_Tqq} into acccount, we can simplify this to
    \begin{align*}
      G_1 T_{q,q'}G_1^T + G_1 T_{q',q} G_1^T&\in\mathrm{Anti}(2),\\
      G_2 T_{q,q'}G_3^T + G_3 T_{q',q} G_2^T&\in\mathrm{Anti}(2),\\
      G_4 T_{q,q'}G_2^T + G_2 T_{q',q} G_4^T&\in\mathrm{Anti}(2).
    \end{align*}
    For any fixed $q,q'$, this is a system of $9$ linear equations over the $18$ coordinates of $T_{q,q'}$ and $T_{q',q}$, which can be seen to be equivalent to
    \begin{equation}\label{eq:orth_cut_Tqqp}
      T_{q,q'}=-T_{q',q}^T.
    \end{equation}

    Combining Equations \eqref{eq:orth_cut_R_Tqq} and \eqref{eq:orth_cut_Tqqp} proves that $R$ is zero and $T$ is antisymmetric, thus $\dot V$ is indeed of the form $VA$ for some $A\in\mathrm{Anti}(p)$, which concludes the proof in the case where $p\equiv 0[3]$.

    In the case where $p\equiv 1[3]$, we define $W_q^{(1)},W_q^{(2)},W_q^{(3)},W_q^{(4)}$ and $X_{q,q'}^{(1)},X_{q,q'}^{(2)},X_{q,q'}^{(3)}$ as previously, for $q\leq \frac{p-4}{3}$, and $q<q'\leq \frac{p-4}{3}$. For $q=\frac{p-1}{3}$, we define six matrices $(W_q^{(i)})_{i=1,\dots,6}$ by
    \begin{equation*}
      W_q^{(i)} = \left( \begin{smallmatrix}
          0_{2\times (p-4)}&H_i
        \end{smallmatrix} \right),
    \end{equation*}
    with
    \begin{gather*}
      H_1 = \left(\begin{smallmatrix}
          1&0&0&0 \\ 0&1&0&0
        \end{smallmatrix}\right),\qquad
      H_2 = \left(\begin{smallmatrix}
          0&1&0&0 \\ 0&0&1&0
        \end{smallmatrix}\right),\qquad
      H_3 = \left(\begin{smallmatrix}
          0&0&1&0 \\ 0&0&0&1
        \end{smallmatrix}\right), \\
      H_4 = \left(\begin{smallmatrix}
          0&0&0&1 \\ 1&0&0&0
        \end{smallmatrix}\right),\qquad
      H_5 = \left(\begin{smallmatrix}
          \frac{1}{\sqrt{2}}&\frac{1}{\sqrt{2}}&0&0 \\ 0&0&\frac{1}{\sqrt{2}}&\frac{1}{\sqrt{2}}
        \end{smallmatrix}\right),\qquad
      H_6 = \left(\begin{smallmatrix}
          0&0&\frac{1}{\sqrt{2}}&-\frac{1}{\sqrt{2}} \\ \frac{3}{5}&\frac{4}{5}&0&0
        \end{smallmatrix}\right).
    \end{gather*}
    And for $q\leq \frac{p-4}{3},q'=\frac{p-1}{3}$, we define the following four matrices:
    \begin{gather*}
      X_{q,q'}^{(1)} = \left( \begin{smallmatrix}
          0_{2\times 3(q-1)}&G_1&0_{2\times p-3q-4}&H_1
        \end{smallmatrix} \right) / \sqrt{2}, \\
      X_{q,q'}^{(2)} = \left( \begin{smallmatrix}
          0_{2\times 3(q-1)}&G_1&0_{2\times p-3q-4}&H_3
        \end{smallmatrix} \right) / \sqrt{2}, \\
      X_{q,q'}^{(3)} = \left( \begin{smallmatrix}
          0_{2\times 3(q-1)}&G_2&0_{2\times p-3q-4}&H_1
        \end{smallmatrix} \right) / \sqrt{2}, \\
      X_{q,q'}^{(4)} = \left( \begin{smallmatrix}
          0_{2\times 3(q-1)}&G_2&0_{2\times p-3q-4}&H_3
        \end{smallmatrix} \right) / \sqrt{2}.
    \end{gather*}
    We define $V$ as before. Establishing that $\mathcal{M}_p$ is $X_0$-minimally secant at $V$ can be done as previously, the only difference being that we have to write $R,T$ as
\begin{equation*}
      R = \begin{pmatrix} R_1&\dots&R_{(p-1)/3}\end{pmatrix},
    \end{equation*}
    with $R_1,\dots,R_{(p-4)/3}\in\mathbb{R}^{2\times 3},R_{(p-1)/3}\in\mathbb{R}^{2\times 4}$, and
    \begin{center}
      \begin{tikzpicture}
        \node at (-3.4,0) {$T=$};
        \node at (0,0) {$\begin{pmatrix}
            T_{1,1}&\dots&T_{1,\frac{p-1}{3}} \\
            \vdots&&\vdots \\
            T_{\frac{p-1}{3},1}&\dots&T_{\frac{p-1}{3},\frac{p-1}{3}}
          \end{pmatrix}$,} ;
        \draw[decorate,decoration={brace}] (-1.7,1) -- (-1,1);
        \node at (-1.35,1.4) {$3$};
        \draw[decorate,decoration={brace}] (0.5,1) -- (1.5,1);
        \node at (1,1.4) {$4$};
        \draw[decorate,decoration={brace}] (-2.2,0.45) -- (-2.2,0.75);
        \node at (-2.5,0.6) {$3$};
        \draw[decorate,decoration={brace}] (-2.2,-0.8) -- (-2.2,-0.4);
        \node at (-2.5,-0.6) {$4$};
      \end{tikzpicture}
    \end{center}

    Finally, in the case where $p\equiv 2[3]$, we define $W_q^{(1)},W_q^{(2)},W_q^{(3)},W_q^{(4)}$ and $X^{(1)}_{q,q'},X^{(2)}_{q,q'},X^{(3)}_{q,q'}$ as before for $q\leq \frac{p-2}{3}$ and $q<q'\leq \frac{p-2}{3}$. For $q=\frac{p+1}{3}$, we define only three matrices $(W_q^{(i)})_{i=1,2,3}$:
    \begin{equation*}
      W_q^{(i)} = \left( \begin{smallmatrix}
          0_{2\times (p-2)}&J_i
        \end{smallmatrix} \right),
    \end{equation*}
    with
    \begin{gather*}
      J_1 = \left(\begin{smallmatrix}
          1&0 \\ 0&1
        \end{smallmatrix}\right),\qquad
      J_2 = \left(\begin{smallmatrix}
          0&1 \\ -1&0
        \end{smallmatrix}\right),\qquad
      J_3 = \left(\begin{smallmatrix}
          \frac{1}{\sqrt{2}}&\frac{1}{\sqrt{2}} \\ \frac{1}{\sqrt{2}}&-\frac{1}{\sqrt{2}}
        \end{smallmatrix}\right).
    \end{gather*}
    For $q\leq \frac{p-2}{3}$, $q'=\frac{p+1}{3}$, we set
    \begin{gather*}
      X_{q,q'}^{(1)} = \left( \begin{smallmatrix}
          0_{2\times 3(q-1)}&G_1&0_{2\times p-3q-2}&J_1
        \end{smallmatrix} \right) / \sqrt{2}, \\
      X_{q,q'}^{(2)} = \left( \begin{smallmatrix}
          0_{2\times 3(q-1)}&G_2&0_{2\times p-3q-2}&J_1
        \end{smallmatrix} \right) / \sqrt{2}.
    \end{gather*}
    We conclude as before.

    The case $d=3$ can be dealt with in the same way as $d=2$, but is even more technical. The easiest thing to do (although maybe not the most elegant one) is to distinguish $12$ cases, depending on the congruency of $p$ modulo $12$. To avoid pages of definitions, we only focus on the case where $p\equiv 0 [12]$ and, even in this case, only provide a sketch of proof.

    For any $q\leq \frac{p}{12}$, we define $19$ matrices $(W_q^{(i)})_{i\leq 19}$ of size $3\times p$, by
    \begin{equation*}
      W_q^{(i)} = \left( \begin{smallmatrix}
          0_{3\times 12(q-1)}&G_i&0_{3\times (p-12q)}
        \end{smallmatrix} \right),
    \end{equation*}
    for matrices $G_1,\dots,G_{19}\in\R^{3\times 12}$ suitably chosen\footnote{By analycity arguments, one can see that a ``generic'' choice of $G_1,\dots,G_{19}$ such that $G_iG_i^T=I_3$ works.}.
    Then, for any $q,q'\leq \frac{p}{12}$, with $q<q'$, one defines $24$ matrices $(X_{q,q'}^{(i)})_{i\leq 24}$
    \begin{gather*}
      X_{q,q'}^{(i)} = \left( \begin{smallmatrix}
          0_{3\times 12(q-1)}&G'_i&0_{3\times 12(q'-q-1)}&G''_i&0_{3\times (p-12q')}
        \end{smallmatrix} \right),
    \end{gather*}
    for appropriate $G'_i,G''_i$.

    If we divide $R$ and $T$ into, respectively, $3\times 12$ and $12\times 12$ blocks, we can check in the exact same way as in the case $d=2$ (using a computer to solve the linear systems) that, if $\dot V = U_0 R + V T$ is in $T_V\mathcal{M}_p$, then $\dot V = VA$ for some $A=\mathrm{Anti}(p)$.

    \subsection{Proof of Corollary \ref{cor:product_of_spheres}\label{app:product_of_spheres}}

    As in the previous two subsections, it suffices to show fix $p$ such that $\frac{p(p+1)}{2}+p\leq S$ and show that the assumptions in Theorem \ref{theorem:main} are verified. The second assumption is elementary. We show the first and third ones by deducing them from the same assumptions in the \textit{MaxCut} case (Subsection \ref{app:maxcut}).

    Let us set, for any $q\in\mathbb{N}$,
    \begin{align*}
      \mathrm{Ins}_q:V\in\mathbb{R}^{S\times q}
      &\to
      \left(\begin{smallmatrix}
          \begin{smallmatrix}
            V_{1,1}&\dots&V_{1,q}
          \end{smallmatrix}\\
          0_{d_1-1,q}\\
          \begin{smallmatrix}
            V_{2,1}&\dots&V_{2,q}
          \end{smallmatrix}\\
          0_{d_2-1,q} \\
          \vdots
        \end{smallmatrix}\right) \in\mathbb{R}^{D\times q}.
      % \mathrm{Del}_q:V\in\mathbb{R}^{D\times q}
      % &\to \left(\begin{smallmatrix}
      %     V_{1,1}&\dots&V_{1,q}\\
      %     V_{d_1+1,1}&\dots&V_{d_1+1,q}\\
      %     V_{d_1+d_2+1,1}&\dots&V_{d_1+d_2+1,q}\\
      %     &\vdots&
      %     \end{smallmatrix}\right) \in\mathbb{R}^{S\times q}.
    \end{align*}

    We set
    \begin{equation*}
      U_0 = \mathrm{Ins}_1(1_{S,1})\quad\mbox{and}\quad
      X_0 = U_0 U_0^T.
    \end{equation*}
    It is an extreme point of $\mathcal{C}$ with rank $1$, so the first hypothesis holds true.

    To establish the third hypothesis, let $V^{(MaxCut)}\in\mathbb{R}^{S\times p}$ be the matrix defined in Subsection \ref{app:maxcut} (with $n=S$). We set
    \begin{equation*}
      V = \mathrm{Ins}_p(V^{(MaxCut)}).
    \end{equation*}
    If $\dot V$ is an element of $T_V\mathcal{M}_p$ such that
    \begin{equation*}
      \mathrm{Range}(\dot V)
      \subset \mathrm{Range}(X_0) + \mathrm{Range}(V),
    \end{equation*}
    we check that $\dot V$ must be of the form $\mathrm{Ins}_p(\dot V^{(MaxCut)})$ for some $\dot V^{(MaxCut)}\in T_V\mathcal{M}^{(MaxCut)}_p$ such that
    \begin{equation*}
      \mathrm{Range}(\dot V^{(MaxCut)})
      \subset \mathrm{Range}(U_0^{(MaxCut)}) + \mathrm{Range}(V^{(MaxCut)}).
    \end{equation*}
    (Here, $\mathcal{M}_p^{(MaxCut)}$ is the feasible set of Problem \eqref{eq:P} in the \textit{MaxCut} case, for $n=S$, and we define $U_0^{(MaxCut)}=1_{S,1}$.)

    From the proof in Subsection \ref{app:maxcut}, there exists an antisymmetric $A\in\mathrm{Anti}(p)$ such that $\dot V^{(MaxCut)}=V^{(MaxCut)}A$. For this $A$, we have
    \begin{equation*}
      \dot V = \mathrm{Ins}_p(\dot V^{(MaxCut)})
      = \mathrm{Ins}_p(V^{(MaxCut)}) A
      = VA.
    \end{equation*}
    Therefore, $\mathcal{M}_p$ is $X_0$-minimally secant at $V$, which concludes the proof.

%\bibliography{../../../biblio/bib_articles.bib,../../../biblio/bib_proceedings.bib,../../../biblio/bib_livres.bib,../../../biblio/bib_misc.bib}

\end{document}